\documentstyle[12pt]{article}

\newcommand{\sect}[1]{\setcounter{equation}{0}\section{#1}}

\newcommand{\eq}{\begin{equation}}
\newcommand{\eqa}{\begin{eqnarray}}
\newcommand{\en}{\end{equation}}
\newcommand{\ena}{\end{eqnarray}}
\newcommand{\enn}{\nonumber \end{equation}}

\def\sk{\vskip .4cm}
\def\noi{\noindent}
\def\om{\omega}
\def\al{\alpha}
\def\be{\beta}
\def\ga{\gamma}
\def\Ga{\Gamma}
\def\del{\delta}

\def\linv{{1 \over \lambda}}
\def\rinv{{1\over {r-r^{-1}}}}

\def\epsi{\varepsilon}

\def\de{\delta}

\def\part{\partial}

\def\R#1#2{ R^{#1}_{~~~#2} }
\def\PA#1#2{ P^{#1}_{A~~#2} }

\def\Rp#1#2{ (R^+)^{#1}_{~~~#2} }

\def\Rm#1#2{ (R^-)^{#1}_{~~~#2} }
\def\Rinv#1#2{ (R^{-1})^{#1}_{~~~#2} }

\def\Rpm#1#2{(R^{\pm})^{#1}_{~~~#2} }

\def\Rh{{\hat R}}

\def\Rhat#1#2{ \Rh^{#1}_{~~~#2} }

\def\Rhatinv#1#2{ (\Rh^{-1})^{#1}_{~~~#2} }

\def\ff#1#2#3{f_{#1~~~#3}^{~#2}}
\def\MM#1#2#3{M^{#1~~~\!\!\!#3}_{~#2}}
\def\cchi#1#2{\chi^{#1}_{~#2}}

\def\U#1#2#3#4#5#6#7#8{U^{~#2~#4}_{#1~#3}|^{#5~#7}_{~#6~#8}}

\def\bu{\bullet}
\def\ci{\circ}

\def\T#1#2{ T^{#1}_{~~#2} }
\def\t#1#2{ T^{#1}_{~~#2} }
\def\caM{{t}}
\def\M#1#2{ t^{#1}_{~\,#2} }

\def\rminus{r^{-1}}

\def\D{\Delta}

\def\Mat#1#2#3#4#5#6#7#8#9{\left( \matrix{
     #1 & #2 & #3 \cr
     #4 & #5 & #6 \cr
     #7 & #8 & #9 \cr
   }\right) }

\def\ep{\epsi^{\prime}}
\def\kp{\kappa^{\prime}}

\def\IGLqrN{IGL_{q,r}(N)}

\def\ISOqrN{ISO_{q,r}(N)}

\def\SqrNt{S_{q,r}(N+2)}

\def\Tc{{\cal T}}

\def\Dtwo{\Delta_{N+2}}
\def\epsitwo{\epsi_{N+2}}
\def\kappatwo{\kappa_{N+2}}

\def\Lpm#1#2{L^{\pm #1}_{~~~#2}}

\def\LLpm{L^{\pm}}

\def\LLp{L^{+}}
\def\LLm{L^{-}}
\def\Lp#1#2{L^{+ #1}_{~~~#2}}
\def\Lm#1#2{L^{- #1}_{~~~#2}}
\def\Lc{{{\cal L}^+}}

\def\P{P}

\def\n2{{{N+1} \over 2}}
\def\ap{a^{\prime}}
\def\bp{b^{\prime}}
\def\cp{c^{\prime}}
\def\dpr{d^{\prime}}

\def\Ntwo{{N\over 2}}

\def\square{{\,\lower0.9pt\vbox{\hrule \hbox{\vrule height 0.2 cm
\hskip 0.2 cm \vrule height 0.2 cm}\hrule}\,}}

\def\SO{SO_{q,r}(N+2)}
\def\ISO{ISO_{q,r}(N)}
\def\U{U_{q,r}(so(N+2))}
\def\IU{U_{q,r}(iso(N))}

\def\ISOqrN{ISO_{q,r}(N)}

\def\SqrNt{S_{q,r}(N+2)}

\def\H{H^\bot}
\def\le{\langle}
\def\re{\rangle}

\def\limrone{\lim_{r \rightarrow 1}}
\def\linv{{1 \over \lambda}}
\def\sma#1{\mbox{\footnotesize #1}}
\def\cvd{{\vskip -0.49cm\rightline{$\Box\!\Box\!\Box$}}\sk}

\def\Re{R}

\begin{document}

\begin{titlepage}
\rightline{LBNL-40330}
\rightline{IFUP-TH 63/95}
\rightline{DFTT-73/95}
\rightline{August 1996}
\vskip 1.2em

\begin{center}
{{\bf Universal Enveloping Algebra and}
\\[1em]
{\bf Differential Calculi on Inhomogeneous Orthogonal
\mbox{\boldmath $q$}-Groups}}
\\[2.8em]
Paolo Aschieri \\[1em]
{\sl Theoretical Physics Group, Physics Division\\
Lawrence Berkeley National Laboratory, 1 Cyclotron Road \\
Berkeley, California 94720, USA.}\\[3em]
Leonardo Castellani\\[1em]
{\sl Dipartimento di Scienze e Tecnologie
avanzate${~}^\diamond$, Universit\`a di Torino}\\[.5 em]
{\sl Dipartimento di Fisica Teorica\\
and\\Istituto Nazionale di
Fisica Nucleare\\
Via P. Giuria 1, 10125 Torino, Italy.}  \\[3em]
\end{center}

\begin{abstract}
\noi ~~We review the construction of the multiparametric quantum  
group
$\ISO$  as a projection from $\SO$ and show that it is a
bicovariant bimodule over $SO_{q,r}(N)$.
The universal enveloping algebra $U_{q,r}(iso(N))$, characterized  as
the Hopf
algebra of regular functionals on $\ISO$, is found
as a Hopf subalgebra of $\U$ and  is shown to be a bicovariant
bimodule over $U_{q,r}(so(N)).$

An $R$-matrix formulation of $U_{q,r}(iso(N))$
is given and  we prove the pairing
$U_{q,r}(iso(N))\leftrightarrow \ISO$.
We analyze  the subspaces of $U_{q,r}(iso(N))$ that define
bicovariant differential calculi on $\ISO$.

\end{abstract}

\vskip .8cm

{}~~~~~~~~~~~~~~~~~~\small{Subj. ind. class. 17B37 81R50
16W30.}~~~~~~~~~~~q-alg/9705023
\vskip .2cm
\noi \hrule
\vskip.2cm
\noi {\small ${}^\diamond$ II Facolt\`a di Scienze Matematiche,  
Fisiche e
Naturali,
sede di Alessandria}

\noi {\small {\sl e-mail}: aschieri@theorm.lbl.gov,
castellani@to.infn.it}

\end{titlepage}
\newpage
\setcounter{page}{1}


\sect{Introduction}

A noncommutative space-time, with a deformed Poincar\'{e} symmetry  
group,
is an interesting geometric background for the study of  Minkowski
space-time
physics and,
in particular, of  Einstein-Cartan gravity theories
\cite{noncom}, \cite{Cas2}.
In this perspective it is natural to investigate
inhomogeneous orthogonal quantum groups,
their quantum Lie algebras and more generally their differential
structure.

In this paper we review the multiparametric $R$-matrix
formulation of $\ISO$ as a
projection from $\SO$ \cite{inson} emphasizing the analogy with the
classical construction. We also show that $\ISO$ is a bicovariant
bimodule over $SO_{q,r}(N)$, freely generated by the translation
elements $x^a$ plus the dilatation element associated to $\ISO$.
We then construct and analize the universal enveloping algebra $\U$,
defined as the algebra of regular functionals \cite{FRT}
on the multiparametric homogeneous orthogonal $q$-groups.
The projection procedure $\SO\rightarrow \ISO$,
initiated in \cite{Cas3} and developed in  
\cite{Cas2,Aschieri1,inson},
is here exploited to
obtain $\IU$ as a particular Hopf subalgebra of $\U$, and prove that  
it is
paired to $\ISO$.
A detailed study of $U_{q,r}(iso(N))$ and an $R$-matrix formulation  
is
given.
In complete analogy with the $\ISO$ case we also prove that
$\IU$ is a bicovariant bimodule over $U_{q,r}(so(N))$
and give a basis of right invariant elements that freely generate  
$\IU$.
The universal
enveloping algebras of the inhomogeneous quantum groups
$IGL_{q,r}(N)$,
first studied with a different approach in \cite{Schlie}, can be
derived in  a similar way.

The quantum Lie algebras of $\ISO$ are subspaces (adjoint submodules)  
of
$\IU$, and in the last section we
examine two of them, obtained as ``projections''  from the quantum  
Lie
algebras of $\SO$.
The two associated bicovariant differential calculi are also briefly
presented.
The first  has $N+2$ generators, and is an interesting candidate for  
a
differential calculus on the quantum orthogonal plane in dimension  
$N$.
The second  is obtained with the parametric restriction $r=1$; in the
classical limit $r=q=1$ it reduces to the differential
calculus on the undeformed $ISO(N)$. This section does not rely on  
the
technical
parts of Section 4 and 5; these may be skipped by the reader  
interested
mainly
in the differential calculi on $ISO_{q,r}(N)$.

In this article, all the
properties of the
quantum inhomogeneous $\ISO$ group, its universal enveloping
algebra and its differential
calculus are derived from the known properties
of the homogeneous ``parent'' structure.
The main logical steps of this derivation are independent
from the $q$-group considered, and the projection
procedure may be applied to
investigate more general quotients of the
$A,B,C,D$ $q$-groups, as for example deformed parabolic groups.

\sect{$SO_{q,r}(N)$ multiparametric quantum group}

The $SO_{q,r}(N)$ multiparametric
quantum group is freely generated by the noncommuting
matrix elements $\T{a}{b}$ (fundamental representation $a,b= 1,\ldots  
N$)
and
the unit element $I$, modulo the relation $\mbox{det}_{q,r}T=I$ and  
the
quadratic $RTT$ and $CTT$ (othogonality)
relations discussed below. The noncommutativity is controlled by the  
$R$
matrix:
\eq
\R{ab}{ef} \T{e}{c} \T{f}{d} = \T{b}{f} \T{a}{e} \R{ef}{cd}
\label{RTT}
\en
which satisfies the quantum Yang-Baxter equation
\eq
\R{a_1b_1}{a_2b_2} \R{a_2c_1}{a_3c_2} \R{b_2c_2}{b_3c_3}=
\R{b_1c_1}{b_2c_2} \R{a_1c_2}{a_2c_3} \R{a_2b_2}{a_3b_3}, \label{QYB}
\en
a sufficient condition for the consistency of the
``$RTT$" relations (\ref{RTT}).  The $R$-matrix components
$\R{ab}{cd}$ depend continuously on a (in general complex)
set of  parameters $q_{ab},r$.  For $q_{ab}= r$ we
recover the uniparametric orthogonal group $SO_r(N)$ of ref.  
\cite{FRT}.
Then
$q_{ab} \rightarrow 1, r \rightarrow 1$ is the classical limit for
which
$\R{ab}{cd} \rightarrow \de^a_c \de^b_d$ : the
matrix entries $\T{a}{b}$ commute and
become the usual entries of the fundamental representation.
The multiparametric $R$ matrices for the $A,B,C,D$ series can be
found in \cite{Schirrmacher}  (other ref.s on multiparametric
$q$-groups are given in \cite{multiparam1,multiparam2}). For
the orthogonal  case they read
(we use the same notations of \cite{inson}):
\eq
\begin{array}{ll}
\R{ab}{cd}=&\delta^a_c \delta^b_d [{r\over q_{ab}} +
(r-1) \delta^{ab}+(\rminus - 1)
\de^{a\bp}] (1-\de^{a n_2})
+\de^a_{n_2} \de^b_{n_2} \de^{n_2}_c \de^{n_2}_d  \\
&+(r-r^{-1})
[\theta^{ab} \delta^b_c \de^a_d -
\theta^{ac} r^{\rho_a - \rho_c}
\de^{\ap b} \de_{\cp d}]
\end{array}
\label{Rmp}
\en
\noi where $\theta^{ab}=1$ for $a> b$
and $\theta^{ab}=0$ for $ a \leq b$; we define
 $n_2 \equiv \n2$ and primed indices as $\ap \equiv N+1-a$.
The terms with the index $n_2$ are present
only in the $B_n$ case: $N=2n+1$. The
$\rho_a$ vector is  given by:
\eq
(\rho_1,...\rho_N)=\left\{ \begin{array}{ll}
         (\Ntwo -1, \Ntwo -2,...,
{1\over 2},0,-{1\over 2},...,-\Ntwo+1)
                   & \mbox{for $B_n[SO(2n+1)]$} \\
           (\Ntwo -1,\Ntwo -2,...,1,0,0,-1,...,-\Ntwo+1) & \mbox{for
$D_n[SO(2n)]$}
                                             \end{array}
                                    \right.
\en
Moreover the following relations reduce the number of independent
$q_{ab}$ parameters \cite{Schirrmacher}:
\eq
q_{aa}=r,~~q_{ba}={r^2 \over q_{ab}}; \label{qab1}
\en
\eq
q_{ab}={r^2 \over q_{a\bp}}={r^2 \over q_{\ap b}}=q_{\ap\bp}
 \label{qab2}
\en
\noi where (\ref{qab2}) also implies $q_{a\ap}=r$. Therefore
 the $q_{ab}$ with $a < b \leq {N\over 2}$ give all the $q$'s.
\sk
It is useful to list the nonzero complex components
of the $R$ matrix (no sum on repeated indices):
\eqa
& &\R{aa}{aa}=r , ~~~~~~~~~~~~~~\mbox{\footnotesize
 $a \not= {n_2}$ } \cr
& &\R{a\ap}{a\ap}=r^{-1} ,  ~~~~~~~~~~\mbox{\footnotesize
 $a \not= {n_2}$ } \cr
& &\R{{n_2}{n_2}}{{n_2}{n_2}}= 1\cr
& &\R{ab}{ab}={r \over q_{ab}} ,~~~~~~~~~~~~\mbox{\footnotesize
 $a \not= b$,  $\ap \not= b$}\cr
& &\R{ab}{ba}=r-r^{-1} , ~~~~~~~\mbox{\footnotesize
$a>b, \ap \not= b $}\cr
& &\R{a\ap}{\ap a}=(r-r^{-1})(1- r^{\rho_a-\rho_{\ap}}) ,
{}~~~~~~\mbox{\footnotesize
$a>\ap $}\cr
& &\R{a\ap}{b \bp}=-(r-r^{-1})
r^{\rho_a-\rho_b} , ~~~~~~~
\mbox{\footnotesize $~~a>b ,~ \ap \not= b $}
\label{Rnonzero}
\ena

{\sl Remark 2.1 :}  The matrix $R$ is upper triangular (i.e.
$\R{ab}{cd}=0$ if  [$\sma{$a=c$}$ and $\sma{$b<d$}$] or
$\sma{$a<c$}$)
and has the following properties:
\eq
R^{-1}_{q,r}=R_{q^{-1},r^{-1}}~~;~~~(R_{\!q,r})^{ab}{}_{cd}=(R_{\!q,r} 
)^{\
cp\dpr}{}_{\ap\bp}~~;~~~
(R_{q,r})^{ab}{}_{cd}=(R_{\!p,r})^{dc}{}_{ba}
\label{Rprop1}
\en
where
$q,r$ denote  the set of parameters $q_{ab},r$, and $p_{ab}\equiv  
q_{ba}$.
\sk
\noi The inverse $R^{-1}$ is defined by
$\Rinv{ab}{cd} \R{cd}{ef}=\de^a_e \de^b_f=\R{ab}{cd}
\Rinv{cd}{ef}$.
The first equation in  (\ref{Rprop1}) implies
that for $|q|=|r|=1$, ${\bar R}=R^{-1}$.
\sk
{\sl Remark 2.2 :}
The characteristic equation and the projector decomposition
of $\Rh_{q,r}$, where  $\Rh^{ab}{}_{cd}\equiv R^{ba}{}_{cd}$
are the same as in
the uniparametric case \cite{multiparam1,inson};
in particular the projectors read:
\eqa
&&\!\!\!\!\!\!\!\!\!\!\!P_S={1 \over {r+\rminus}}
[\Rh+\rminus I-(\rminus+r^{1-N})P_0]\,,~~\!
P_A={1 \over {r+\rminus}} [-\Rh+rI-(r- r^{1-N})P_0]\,,\nonumber\\
&&\!\!\!\!\!\!\!\!\!\!\!P_0= (C_{ab} C^{ab})^{-1}K~,~~\mbox{ where }~
K^{ab}_{~~cd}\equiv C^{ab} C_{cd}~.
\label{proiett}
\ena
\sk
Orthogonality  conditions are
imposed on the elements $\T{a}{b}$, consistently
with  the $RTT$ relations (\ref{RTT}):
\eq
C^{bc} \T{a}{b}  \T{d}{c}= C^{ad} I~~~;~~~~
C_{ac} \T{a}{b}  \T{c}{d}=C_{bd} I \label{Torthogonality}
\en
\noi where the (antidiagonal) metric is :
\eq
C_{ab}= r^{-\rho_a} \de_{a\bp} \label{metric}
\en
\noi and its inverse $C^{ab}$
satisfies $C^{ab} C_{bc}=\de^a_c=C_{cb} C^{ba}$.
We see
that the matrix elements of the metric
and the inverse metric coincide: $C^{ab}=C_{ab}$;
notice also the symmetry $C_{ab}=C_{\bp\ap}$.

The consistency of (\ref{Torthogonality}) with the $RTT$ relations
is due to the identities:
\eq
C_{ab} \Rhat{bc}{de} = \Rhatinv{cf}{ad} C_{fe} \label{crc1}
\en
\eq
 \Rhat{bc}{de} C^{ea}=C^{bf} \Rhatinv{ca}{fd} \label{crc2}
\en
\noi These identities
 hold also for $\Rh \rightarrow \Rh^{-1}$ and can be proved using
the explicit expression (\ref{Rnonzero}) of $R$.
We also note the useful relations
\eq
C_{ab}\Rhat{ab}{cd}=r^{1-N}C_{cd} ,~~~
C^{cd}\Rhat{ab}{cd}=r^{1-N}C^{ab}  \label{CR}~,
\en
and
\eq
 \R{ab}{cc'}=C^{ab}C_{cc'}\;,~
\R{aa'}{cd}=C^{aa'}C_{cd} ~~\mbox{ for } ~\sma{$a>c$}~.
\label{foraa'}
\en
The co-structures of the orthogonal multiparametric quantum
group have the same form as in the uniparametric case:
the coproduct
$\D$, the counit $\epsi$ and the coinverse $\kappa$ are given by
\eqa
& & \D(\T{a}{b})=\T{a}{b} \otimes \T{b}{c}  \label{cos1} \\
& & \epsi (\T{a}{b})=\delta^a_b\\
& & \kappa(\T{a}{b})=C^{ac} \T{d}{c} C_{db}
\label{cos2}
\ena

In order to define the quantum determinant $\mbox{det}_{q,r}T$
we introduce the
orthogonal $N$-dimensional quantum plane of coordinates $x^a$
that satisfy the $q$-commutation relations  
${P_A}^{ab}_{~cd}x^cx^d=0$.
We then consider the algebra of exterior forms
$dx^1,dx^2,\ldots dx^N$ defined by: ${P_S}^{ab}_{~cd}dx^cdx^d=0$ and
${P_0}^{ab}_{~cd}dx^cdx^d=0$ i.e. [use (\ref{proiett})]:
$dx^adx^b=-rR^{ba}_{~cd}dx^cdx^d$. There is a natural action $\de$
of the orthogonal quantum group  on the exterior algebra
(that becomes a left comodule):
\[\de(dx^a)=\T{a}{c}\otimes dx^c~;~~\de(dx^adx^b ... dx^c)=
\T{a}{d}\T{b}{e} ...\T{c}{f}\otimes
dx^ddx^e ...dx^f~.\]
Generalizing the results of \cite{Fiore}
to the multiparametric case, we find
that any $N$-dimensional form is proportional
to the volume form $dV\equiv dx^1\ldots dx^N$, so that the  
determinant is
uniquely defined by:
\eq
\de(dV)\equiv \mbox{det}_{q,r}T\otimes dV~.
\en
Using
(\ref{Torthogonality})
as in \cite{Fiore} it is immediate to prove that
$(\mbox{det}_{q,r}T)^2=I$; moreover
$\mbox{det}_{q,r}T$ is central
and satisfies
$\D (\mbox{det}_{q,r}T)=
\mbox{det}_{q,r}T\otimes \mbox{det}_{q,r}T$.

To obtain the special orthogonal quantum group $SO_{q,r}(N)$ we  
impose
also the
relation   $\mbox{det}_{q,r}T=I$.
\sk

{\sl Remark 2.3 :}
Using formula (\ref{Rmp}) or (\ref{Rnonzero}),
we find that the $\R{AB}{CD}$  matrix for the
$SO_{q,r}(N+2)$  quantum group can be decomposed
in terms of  $SO_{q,r}(N)$ quantities
as follows (splitting the index {\small A} as
{\small A}=$(\circ, a, \bullet)$, with $a=1,...N$):
\eq
\R{AB}{CD}=\left(  \begin{array}{cccccccccc}
   {}&\circ\circ&\circ\bullet&\bullet
          \circ&\bullet\bullet&\circ d&\bullet d
      &c \circ&c\bullet&cd\\
   \circ\circ&r&0&0&0&0&0&0&0&0\\
   \circ\bullet&0&r^{-1}&0&0&0&0&0&0&0\\
   \bullet\circ&0&f(r)&r^{-1}&0&0&0&0&0&- C_{cd} \lambda
r^{-\rho}\\
\bullet\bullet&0&0&0&r&0&0&0&0&0\\
\circ b&0&0&0&0&{r\over q_{\circ b}} \de^b_d&0&0&0&0\\
\bullet b&0&0&0&0&0&{r\over q_{\bullet b} }
\de^b_d&0&\lambda\de^b_c&0\\
a\circ&0&0&0&0&\lambda\de^a_d&0&{r \over q_{a \circ} } \de^a_c&0&0\\
a\bullet&0&0&0&0&0&0&0&{r\over q_{a \bullet}} \de^a_c&0\\
ab&0&-C^{ba} \lambda r^{-\rho}
&0&0&0&0&0&0&\R{ab}{cd}\\
\end{array} \right) \label{Rbig}
\en
\noi where $\R{ab}{cd}$ is the $R$ matrix for  $SO_{q,r}(N)$,  
$C_{ab}$ is
the corresponding
metric,  $\lambda \equiv r-r^{-1}$,
$\rho={{N}\over 2}~(r^{\rho}=C_{\bullet \circ})$
and $f(r) \equiv \lambda (1- r^{-2\rho})$.

\sk
\sect{$\ISO$ as a projection from $\SO$}
Classically the orthogonal group $SO(N+2)$
is defined  as the set of all linear transformations with unit  
determinant
which preserve
the quadratic form
$(z^0)^2+(z^1)^2+... (z^{N+1})^2$
or equivalently, since we are in the complex plane,
the quadratic form $z^0z^{N+1}+z^1z^{N}+...z^{N+1}z^0$ (use the
transformation $z^A\rightarrow (z^A+iz^{A'})/\sqrt{2} \mbox{ for }
\sma{$A\leq
N/2\;;$}$
$~z^A\rightarrow (z^{A'} -iz^{A})/\sqrt{2} \mbox{ for } \sma{$A >
N/2$}\,;$
$~z^A$ unchanged for  $\sma{$A=A'$}$).
The associated metric is therefore $C_{AB}=\del_{AB'}$ where
\sma{$A,B=0,1,...N+1$} and \sma{$B'\equiv N+1-B$}.

We consider the $ISO(N)$ subgroup of $SO(N+2)$ defined as follows.
Select the subset of matrices in $SO(N+2)$ whose components $T^A{}_B$
read:
\eq
\T{a}{\circ}=\T{\bullet}{b}=\T{\bullet}{\circ}=0~. \label{Tprojected}
\en
The product of two such $SO(N+2)$ matrices gives a $SO(N+2)$ matrix  
with
the same
structure:
\eq
\Mat{\T{\circ}{\circ}}{y}{z}{0}{T}{x}{0}{0}
{{\T{\bullet}{\bullet}}}\cdot
\Mat{T'{}^{\circ}{}_{\circ}}{y'}{z'}{0}{ T'}{x'}{0}
{0}{{T'{}^{\bullet}{}_{\bullet}}}
=\Mat{\T{\circ}{\circ}T'{}^{\circ}{}_{\circ}}{y''}{z''}{0}
{ T\cdot T'}{x''}{0}{0}
{{\T{\bullet}{\bullet}T'{}^{\bullet}{}_{\bullet}}}\label{Matprototto}
\en
where $x^c\equiv T^c{}_{\bullet}$, $y_a\equiv \T{\circ}{a} , z\equiv
\T{\circ}{\bullet} ,
x''=xT'{}^{\bullet}{}_{\bullet}+Tx'$
and $y''=\T{\circ}{\circ}y'+yT'$.  These matrices form a subgroup
of $SO(N+2)$. If we further set  $\T{\circ}{\circ}=\T{\bu}{\bu}=1$  
this
subgroup becomes $ISO(N)$.

The conditions (\ref{Tprojected})  and $\T{A}{B}\in SO(N+2)$
(i.e.
$\T{A}{B}C_{AC}\T{C}{D}=C_{BD}~,~\mbox{det}\T{A}{B}=1$) are  
equivalent to:
\eqa
& & \T{a}{\circ}=\T{\bullet}{b}=\T{\bullet}{\circ}=0~,  
\label{Tproj}\\
& &  
\T{a}{b}C_{ac}\T{c}{d}=C_{bd}~,~\mbox{det}\T{a}{b}=1,\label{auno}\\
& &\T{\circ}{b}=- \T{a}{b} C_{ac} \T{c}{\bullet} \T{\circ}{\circ}~,
{}~\T{\circ}{\bullet}=-{1\over 2}
\T{b}{\bullet}  C_{ba}  \T{a}{\bullet}\T{\circ}{\circ}~,
{}~\T{\circ}{\circ} =(\T{\bullet}{\bullet})^{-1}.\label{uno}
\ena
As expected, there are no constraints on
$x^c\equiv T^c{}_{\bullet}$.
\sk
{\sl{Remark :}} Classically there is an easier way to recover
$ISO(N)$, namely starting from $SO(N+1)$.
At the quantum level the $R$-matrix of $SO_{q,r}(N)$ is only  
contained in
$\SO$,
see {\sl Remark 2.3}. This explains why we have considered
this bigger group.
\sk
Since $ISO(N)$ is a subgroup of $SO(N+2)$ the algebra
$Fun(ISO(N) )$ of regular
functions from $ISO(N)$ to {\boldmath $C$}
 will be obtained from
$Fun(SO(N+2))$ as a quotient, whose canonical projection we name $P$.
Let us now consider the elements $T^A{}_B$
as {\sl functions} on the $SO(N+2)$ group manifold: they define
the fundamental
representation of $SO(N+2)$.
Since $\forall g\in ISO(N)~~$,
$\T{a}{\circ}(g)=\T{\bullet}{b}(g)=\T{\bullet}{\circ}(g)=0~,$
we can write
\eq
Fun(ISO(N))={Fun(SO(N+2))\over H}\label{quoziente}
\en
where $Fun(SO(N+2))$ is generated by $T^A{}_B$ and $H$ is the left
and right ideal generated by the functions
$\T{a}{\circ}~;~\T{\bullet}{b}~;~\T{\bullet}{\circ}~.$
Therefore  $Fun(ISO(N))$ is
generated by the functions $P(T^A{}_B)$ where $P$ is the canonical
projection associated to $H~:~~P(\T{a}{\circ})=P(\T{\bullet}{b})=
P(\T{\bullet}{\circ})=0$; more explicitly it is
generated by the elements
$T^{A}{}_{B}$
modulo the relations (\ref{Tproj})-(\ref{uno}).
\sk

The above construction can be carried over to the quantum group  
level.
In this case the elements $T^A{}_B$ are abstract generators of  
$\SO\equiv
Fun_{q,r}(SO(N+2))$
and we have $\ISO\equiv Fun_{q,r}(ISO(N))=\SO/H$ because the
ideal $H$ is a Hopf ideal i.e.

i) $H$ is a two-sided ideal in $\SqrNt$,

ii) $H$ is a co-ideal, i.e.
\eq
\Dtwo (H) \subseteq H \otimes \SO + \SO \otimes H;~~\epsitwo
(H)=0
\label{coideal}
\en

iii) $H$ is compatible with $\kappatwo$:
\eq
\kappatwo (H)\subseteq H \label{Hideal}
\en
where the suffix \sma{$N+2$} refers to the costructures of $\SO$.
It should be clear that $\ISO$ is {\sl not} a subalgebra,
nor a Hopf subalgebra of $\SO$; that is why we distinguish
with a suffix between the costructures of $\ISO$ and of $\SO$.

Following \cite{inson}
the projection $P~:~\SO\rightarrow\SO/H$ is a Hopf algebra  
epimorphism,
and defining the projected matrix elements $\M{A}{B}=P(\T{A}{B})$,  
where
$\T{A}{B}$ are the
$\SO$ generators, we have the:
\sk
{\sl Theorem 3.1} The quantum group $\ISO$ is generated by the matrix
entries
\eq
{\caM}
\equiv\Mat{P(\T{\circ}{\circ})}{P(y)}{P(z)}{0}{{P(\T{a}{b})}}{P(x)}{0} 
{0}{P(\T{\bullet}{\bullet})}
{}~~\mbox{ and the unity }~I
\en

modulo the ``$R\caM\caM$''  and ``$C\caM\caM$'' relations
\eq
\R{AB}{EF} \caM^{E}{}_{C} \caM^{F}{}_{D} = \caM^{B}{}_{F}  
\caM^{A}{}_{E}
\R{EF}{CD}~,
\label{RTTISO}
\en
\eq
C^{BC} \caM^{A}{}_{B}  \caM^{D}{}_{C}= C^{AD} ~~;~~ C_{AC}  
\caM^{A}{}_{B}
\caM^{C}{}_{D}=C_{BD} \label{CMMbig}~,
\en
where $R$ and $C$ are the multiparametric $R$-matrix and metric of
$\SO$, respectively.
The co-structures are the same as in (\ref{cos1})-(\ref{cos2}), with
$\t{A}{B} $ instead of $\T{a}{b}$. \cvd

Relations (\ref{RTTISO}) and (\ref{CMMbig}) explicitly read:
\eqa
& &\R{ab}{ef} \T{e}{c} \T{f}{d} = \T{b}{f} \T{a}{e} \R{ef}{cd}
\label{PRTT11}\\
& &\T{a}{b} C^{bc} \T{d}{c}=C^{ad} I \label{PRTT31}\\
& &\T{a}{b} C_{ac} \T{c}{d} = C_{bd} I \label{PRTT32}
\ena
\eqa
& &\T{b}{d} x^a={r \over q_{d\bullet}} \R{ab}{ef} x^e \T{f}{d}
\label{PRTT33}\\
& &\PA{ab}{cd} x^c x^d=0 \label{PRTT13}\\
& &\T{b}{d} v={q_{b\bullet}\over q_{d\bullet}} v \T{b}{d}\\
& &x^b v=q_{b \bullet} v x^b \label{PRTT15}\\
& & uv=vu=I \label{PRTT21}\\
& &u x^b=q_{b\bullet} x^bu \label{PRTT22}\\
& &u \T{b}{d}={q_{b\bullet}\over q_{d\bullet}} \T{b}{d} u
\label{PRTT24}
\ena
\eq
y_b=-r^{\rho} \T{a}{b} C_{ac} x^c u \label{ipsilon}
\en
\eq
z=-{1\over {(r^{-{N\over 2}}+r^{{N\over 2}-2})}} x^b C_{ba} x^a u
\label{PRTT44}
\en
\noi where we have set $P(\T{\ci}{\ci})=u , P(\T{\bu}{\bu})=v$ and,
with abuse of notations, $\T{a}{b}=P(\T{a}{b}), ~x=P(x)$
$y=P(y)$,  $z=P(z)$, and where
$q_{a\bullet}$ are $N$  complex parameters
related by $q_{a\bullet} = r^2 /q_{\ap\bullet}\;,$ with $\ap =  
N+1-a$.
The matrix $P_A$ in eq. (\ref{PRTT13}) is the $q$-antisymmetrizer for  
the
orthogonal quantum group, see (\ref{proiett}).
The last two relations (\ref{ipsilon}) - (\ref{PRTT44})
are constraints, showing that
the $t^{A}{}_{B}$ matrix elements
are really a {\sl redundant} set. This redundance
is necessary if we want an $R$-matrix formulation giving the
$q$-commuations
of the $\ISOqrN$  generators.
  Remark that,
in the $R$-matrix formulation
for $\IGLqrN$, {\sl all}  the $\M{A}{B}$
are independent \cite{Cas3,Aschieri1}. Here
we can take as independent generators the
elements
\eq
{}~~\T{a}{b} , x^a , v , u\equiv v^{-1}
\mbox{ and the identity }  I~~~~~(a=1,...N)~.\label{(3.24)}
\en
\indent In the commutative limit $q\rightarrow 1 , r\rightarrow 1$ we
recover the algebra $Fun(ISO(N))$ described in (\ref{quoziente}).
\sk
{\sl Note 3.1 :} From the commutations
(\ref{PRTT22}) - (\ref{PRTT24})
 we see that
one can set $u=I$ only when $q_{a\bullet}=1$ for all $a$.
{}From $q_{a\bullet} = r^2 /q_{\ap\bullet}$, cf. eq. (\ref{qab2}),
this implies also $r=1$.

{\sl Note 3.2 :} eq.s (\ref{PRTT13}) are the multiparametric
orthogonal quantum plane commutations. They follow
from the   $({}^a{}_{\bullet} {}^b{}_{\bullet})$
$Rtt$ components
and (\ref{PRTT44}).

{\sl Note 3.3 : }
The $u, v=u^{-1}$ and $x^a$ elements generate a subalgebra of
$\ISO$ because their commutation relations
do not involve the $\T{a}{b}$ elements.
Moreover these elements can be ordered using
(\ref{PRTT13}) and (\ref{PRTT22}), and the Poincar\'e series of this
subalgebra is the same as that
of the commutative algebra in the $N+1$
symbols $u$, $x^a$ \cite{FRT}.
A linear basis
of this subalgebra is  therefore given by the ordered monomials:
$\zeta^i=u^{i_{\ci}} (x^1)^{i_{1}}...\,
(x^N)^{i_{N}}$
with $i_{\ci}
\in{\bf \mbox{\boldmath$Z$}}$,
$i_1,... i_N\in
{\bf \mbox{\boldmath$N$}}\cup\{0\}$.
Then, using (\ref{PRTT33}) and (\ref{PRTT24}),
a generic element of $\ISO$
can be written as $\zeta^ia_i$ where $a_i\in SO_{q,r}(N)$ and
we conclude that $\ISO$ is a right $SO_{q,r}(N)$--module
generated by the ordered
monomials $\zeta^i.$

One can show that as in the classical case the expressions
$\zeta^ia_i$
are unique:
$\zeta^ia_i=0 \Rightarrow a_i=0 \;\forall\;\sma{$i$}$, i.e.
that $\ISO$ is a {\sl free} right $SO_{q,r}(N)$--module; moreover
(at least when  $q_{a\bu}=r\;\forall a$)
$\ISO$ is a bicovariant bimodule over $SO_{q,r}(N)$.
The proofs of these statements follow the same steps as those
given after {\sl Note 5.4}, and are left to the reader.
Similarly, writing  $a_i\zeta^i$
instead of $\zeta^ia_i$, we have that $\ISO$ is the free left
$SO_{q,r}(N)$--module generated by the
$\zeta^i$.
\sk


\sect{Universal enveloping algebra $\U$}

We construct the universal enveloping algebra $\U$ of $\SO$ as  the
algebra of regular functionals \cite{FRT} on $\SO$.

It is the algebra over $\mbox{\boldmath $C$}$ generated
by the counit $\epsi$ and by the functionals $\LLpm $
defined
by their value on the matrix elements $\T{A}{B}$  :
\eq
\Lpm{A}{B} (\T{C}{D})= \Rpm{AC}{BD}, \label{LonT}
\en
\eq
\Lpm{A}{B} (I)=\de^A_B \label{LonI}
\en
\noi with
\eq
\Rp{AC}{BD} \equiv \R{CA}{DB} \label{Rplus}~~;~~~
\Rm{AC}{BD} \equiv \Rinv{AC}{BD}~.
\en
To extend the definition (\ref{LonT})
to the whole algebra $\SO$ we set
\eq
\Lpm{A}{B} (ab)=\Lpm{A}{C} (a) \Lpm{C}{B} (b)
{}~~~\forall a,b\in  \SO~.
\label{Lab}
\en
{}From (\ref{LonT}),
using the upper and lower
triangularity of $R^+$ and $R^-$, we see that
$L^+$ is upper
triangular and $L^-$ is lower triangular.
\sk
The commutations between $\Lpm{A}{B}$
and $\Lpm{C}{D}$ are induced by (\ref{QYB}) :
\eq
R_{12} \LLpm_2 \LLpm_1=\LLpm_1 \LLpm_2 R_{12} \label{RLL}~,
\en
\eq
R_{12} \LLp_2 \LLm_1=\LLm_1 \LLp_2 R_{12}~, \label{RLpLm}
\en
\noi where as usual the product $\LLpm_2 \LLpm_1$
is the convolution
product $\LLpm_2 \LLpm_1 \equiv (\LLpm_2 \otimes \LLpm_1)\D$.
\sk
The $\Lpm{A}{B}$ elements satisfy orthogonality conditions
analogous to (\ref{Torthogonality}):
\eqa
& &C^{AB} \Lpm{C}{B} \Lpm{D}{A} = C^{DC} \epsi\label{CLL1}\\
& &C_{AB} \Lpm{B}{C} \Lpm{A}{D} = C_{DC} \epsi \label{CLL2}
\ena
\noi as can be verified by applying them to the $q$-group
generators  and using (\ref{crc1}), (\ref{crc2}).
They provide the inverse for the matrix $\LLpm$:
\eq
[(\LLpm)^{-1}]^A{}_{\!B}=C^{DA} \Lpm{C}{D} C_{BC}
\label{Linverse}
\en
\sk
The co-structures of the algebra generated by the
functionals $L^{\pm}$ and $\epsi$
are defined by the duality (\ref{Lab}):
\eq
\D '(\Lpm{A}{B})(a \otimes b) \equiv \Lpm{A}{B}
(ab)=\Lpm{A}{G}(a) \Lpm{G}{B} (b)
\en
\eq
\ep (\Lpm{A}{B})\equiv \Lpm{A}{B} (I)
\en
\eq
\kp(\Lpm{A}{B})(a)\equiv \Lpm{A}{B} (\kappa (a))
\en
\noi so that
\eqa
\!\!\!\!& & \D ' (\Lpm{A}{B})=\Lpm{A}{G} \otimes
\Lpm{G}{B}\label{copLpm}\\
\!\!\!\!& & \epsi ' (\Lpm{A}{B})=\de^A_B \label{couLpm}\\
\!\!\!\!& &\kp (\Lpm{A}{B})
= [(\LLpm)^{-1}]^A{}_{\!B}
= C^{DA} \Lpm{C}{D} C_{BC}
\label{coiLpm}
\ena

{}From (\ref{coiLpm}) we have that $\kp$ is an inner operation
in the algebra generated by the
functionals $\Lpm{A}{B}$ and $\epsi$, it is then easy to see that  
these
elements generate a Hopf algebra, the Hopf algebra $\U$  of  regular
functionals on the quantum
group $\SO$.
\sk
{\sl Note 4.1 :}
{}From the $CLL$ relations
$\kappa'(\Lpm{A}{B})\Lpm{B}{C}=\Lpm{A}{B}\kappa'(\Lpm{B}{C})=\delta^A_ 
C\ep
si$
we have, using upper-lower triangularity of $\LLpm$:
\eq
\Lpm{A}{A}\kappa'({\Lpm{A}{A}})=
\kappa'(\Lpm{A}{A})\Lpm{A}{A}=\epsi~~\mbox{ i.e. }~~
\Lpm{A}{A}\Lpm{A'}{A'}=\Lpm{A'}{A'}\Lpm{A}{A}=\epsi
\en
As a consequence  $\mbox{det}\LLpm\equiv \Lpm{\circ}{\circ}\Lpm{1}{1}
\Lpm{2}{2}\ldots\Lpm{N}{N}\Lpm{\bullet}{\bullet}
=\epsi$. In the $B_n$ case we also have $\Lpm{n_2}{n_2}=\epsi$.
\sk
{\sl Note 4.2 :}
The $RLL$ relations imply that the subalgebra $U^0$ generated by the
elements
$\Lpm{A}{A}$ and $\epsi$ is commutative (use upper triangularity of  
$R$).
Moreover, from (\ref{copLpm}) the invertible elements $\Lpm{A}{A}$  
are
also group like, and we conclude that
$U^0$ is the group Hopf algebra of the abelian group generated by
$\Lpm{A}{A}$
and $\epsi ~.$
In the classical limit $U^0$ is a maximal commutative subgroup of
$SO(N+2)$.
\sk
\indent {\sl Note 4.3 :}
When $q_{AB}=r$,
the multiparametric
$R$-matrix goes into the uniparametric $R$-matrix and we recover
the standard uniparametric orthogonal quantum groups.
Then
the $\LLpm$ functionals satisfy the further relation:
\eq
\forall\, \sma{A}\;, ~~~~~~
\Lp{A}{A}\Lm{A}{A}
=\epsi \label{epsiaepsi}~,
\en
indeed $\Lp{A}{A}\Lm{A}{A}(a)=
\epsi(a)$ as can be easily seen when $a=\T{A}{B}$
and generalized to any $a\in \SO$ using (\ref{Lab}).
In this case \cite{FRT}
we can avoid to realize the Hopf algebra
$U_r(so(N+2))$ as functionals on  $SO_r(N+2)$ and we can define it
abstractly as
the Hopf algebra generated by  the {\sl symbols} $\LLpm$
and the unit $\epsi$
modulo relations
(\ref{RLL}),(\ref{RLpLm}),(\ref{CLL1}),(\ref{CLL2}), and
(\ref{epsiaepsi}).\sk

As discussed in \cite{FRT} in the uniparametric case, the Hopf  
algebra
$U_r(so(N+2))$
of regular functionals is a Hopf subalgebra of the
orthogonal Drinfeld-Jimbo universal enveloping algebra $U_h$, where
$r=e^h$.
In the general multiparametric case, relation (\ref{epsiaepsi})
does not hold any more. Here we discuss the generalization
of (\ref{epsiaepsi}) and the relation between $\U$ and the  
multiparametric
orthogonal Drinfeld-Jimbo universal enveloping algebra  
$U_h^{(\cal{F})}$.
This latter is the quasitriangular Hopf algebra $U_h^{(\cal{F})}=
(U_h,\Delta^{({\cal{F}})},S,{\cal{R^{(F)}}})$ paired to the
multiparametric orthogonal $q$-group $\SO$.
It is obtained from
$U_h=(U_h,\D,S,{\cal{R}})$ via a twist \cite{multiparam1}.
$U_h^{(\cal{F})}$ has the same algebra structure of $U_h$
(and the same antipode $S$),
while the coproduct $\Delta^{{\cal{(F)}}}$ and the universal element
${\cal{R^{(F)}}}$
belonging to (a completion of) $U_h\otimes U_h$
are determined by the twisting element ${\cal{F}}$
that belongs to (a completion of) a maximal commutative subalgebra of
$U_h\otimes U_h$. We have
\eq
\forall\,\phi\in U_h\,, ~~\D^{{\cal{(F)}}}(\phi)=
{\cal{F}}\D(\phi){\cal{F}}^{-1}~;~~{\cal{R^{(F)}}}=
{\cal{F}}_{21}{\cal{R}}{\cal{F}}^{-1}~;~~
{\cal{R^{(F)}}}(T\otimes T)=R_{q,r}~.
\en
The element ${\cal{F}}$ satisfies:
$(\D^{{\cal{(F)}}}\otimes  
id){\cal{F}}={\cal{F}}_{13}{\cal{F}}_{23}\,,~
(id  
\otimes\D^{{\cal{(F)}}}){\cal{F}}={\cal{F}}_{13}{\cal{F}}_{12}\,,~
{\cal{F}}_{12}{\cal{F}}_{21}=I\,,~
{\cal{F}}_{12}{\cal{F}}_{13}{\cal{F}}_{23}=
{\cal{F}}_{23}{\cal{F}}_{13}{\cal{F}}_{12}\,,$
$\, (\epsi\otimes id){\cal{F}}=(id\otimes \epsi){\cal{F}}
=\epsi\,,~(S \otimes id){\cal{F}}
=(id \otimes S ){\cal{F}}=
{\cal{F}}^{-1},\, \cdot(id\otimes S){\cal{F}}=
\cdot(S\otimes id){\cal{F}} =
\cdot(id \otimes id){\cal{F}} =\epsi\,$;
we explicitly have
\eq
{\cal{F}}(\T{A}{B}\otimes\T{C}{D})= F^{AC}_{~BD}
\en
where $F^{AC}_{~BD}$ is the diagonal matrix
\eq
F=diag (\sqrt{{q_{11}\over r}} ,
\sqrt{{ q_{12}\over r }} , ... ~ \sqrt{{ q_{NN} \over r}})
\en
It is easy to see that the definition of the $\LLpm$ functionals
given in the beginning of this section is equivalent to
the following one:
$\Lp{A}{B}(a)={\cal{R^{(F)}}} ( {a} \otimes \T{A}{B})$ and
 $\Lm{A}{B}(a)={\cal{R^{(F)}}}^{-1}(\T{A}{B}
\otimes {a})$. From
$(\D^{{\cal{(F)}}}\otimes  
id){\cal{R}}={\cal{R}}_{13}{\cal{R}}_{23}\,,~
(id  
\otimes\D^{{\cal{(F)}}}){\cal{R}}={\cal{R}}_{13}{\cal{R}}_{12}\,,~
$ we have $\D^{{\cal{(F)}}}(\Lpm{A}{B})=\Lpm{A}{C} \otimes  
\Lpm{C}{B}$
 and therefore $\D^{{\cal{(F)}}}=\D '$ on $\U$. {}From $(id\otimes S)
({{\cal{R}}})=
(S\otimes id)({{\cal{R}}})={\cal{R}}^{-1}$
it is also easy to see that $S=\kp$ on $\U$
and we conclude that the algebra of regular functionals $\U$ is a
realization [in terms of  functionals on $\SO$] of a
Hopf subalgebra of $U_h^{(\cal{F})}$ with $r=e^h$.
The generalization of (\ref{epsiaepsi}) lies in ${U_h^{{\cal{(F)}}}}$
and not in $\U$, and it is given by
\eq
\forall\, \sma{$A$}~~~~~~~ \Lp{A}{A}\Lm{A}{A}=f_i(\T{A}{A})f^i~~  
\mbox{
where }
{\cal{F}}^4=f_i\otimes f^i~.\label{1effe1}
\en
This relation holds with $\LLpm$ considered as
abstract symbols. It can easily be checked
when $\LLpm$ are realized as functionals:
indeed $\Lp{A}{A}\Lm{A}{A}(a)=
{\cal{F}}^4(\T{A}{A}\otimes a)$ as can be seen when $a=\T{A}{B}$
[use ${\cal{F}}^2( \T{A}{A}\otimes b)={\cal{F}}(\T{A}{A}\otimes b_1)
{\cal{F}}(\T{A}{A}\otimes b_2)$] and generalized to any $a\in \SO$
using ${\cal{F}}(\T{A}{A}\otimes ab)={\cal{F}}(\T{A}{A}\otimes a)
{\cal{F}}(\T{A}{A}\otimes b)$.

In order to characterize the relation between the Hopf algebra
of regular functionals $\U$ and $U_h^{\cal{(F)}}$,
following \cite{FRT}, we extend the group Hopf algebra $U^0$
described in {\sl Note 4.2}
to $\hat{U}^0$
by means of the elements
\footnote{
In the classical limit ${\ell^{\pm}}^A{}_A$ are the tangent vectors
to a maximal commutative subgroup of $SO(N+2)$. They generate a
Cartan subalgebra of the Lie algebra $so(N+2)$.}
${\ell^{\pm}}^A{}_A=\ln \Lpm{A}{A}$. Otherwise stated this means that
in $\hat{U}^0$ we can write $\Lpm{A}{A}=
\mbox{exp}({\ell^{\pm}}^A{}_A)$ where ${\ell^{\pm}}^A{}_A\in  
\hat{U}^0$.
[Explicitly ${\ell^{\pm}}^A{}_A(\T{C}{D})=\ln
({R^{\pm}}^{AC}_{~AC})\,\delta^C_D$,
 ${\ell^{\pm}}^A{}_A(I)=0 $,
${\ell^{\pm}}^A{}_A( ab)={\ell^{\pm}}^A{}_A( a) \epsi (b) + \epsi(a)
{\ell^{\pm}}^A{}_A( b)$
and  $\kp({\ell^{\pm}}^A{}_A)=-{\ell^{\pm}}^A{}_A$ ].
It then follows that ${\cal{F}}$ belongs to (a completion of)
$\hat{U}^0\otimes\hat{U}^0$.
The corresponding extension
$\hat{U}_{q,r} (so(N+2))$ of $\U$, defined as the Hopf algebra  
generated
by the {\sl symbols} $\LLpm$ and $\ell^{\pm}$ modulo relations
(\ref{RLL})-(\ref{CLL2}) and (\ref{1effe1}), is isomorphic -- when  
$r=e^h$
--
to $U_h^{(\cal{F})}\;:\,~\hat{U}_{q,r}(so(N+2))
\cong U_h^{(\cal{F})}$. This relation holds because it
is the twisted version of the
known uniparametric analogue
$\hat{U}_{r}(so(N+2))
\cong U_h$ \cite{FRT,Frenkel}.
\sk

The elements $\LLpm$ [or ${1\over{r-r^{-1}}}(\Lpm{A}{B}-
\de^A_B\epsi$)]
may be seen as the quantum analogue of the tangent vectors; then
the $RLL$ relations are the quantum analogue of the Lie algebra
relations, and we can use the orthogonal $CLL$ conditions
to reduce the number of the $\LLpm$ generators to $(N+2)(N+1)/2$,  
i.e. the
dimension of the classical group manifold.

This we proceed to do; we next  study the $R\LLpm\LLpm$
commutation relations restricted to these $(N+2)(N+1)/2$ generators
and find a set of ordered monomials in the reduced $\LLpm$ that  
linearly
span all
$\hat{U}_{q,r}(so(N+2))$.

We first observe that the commutative subalgebra $\hat{U}^0$
is generated by ${(N+2)/2}$ elements ($N$ even, $N=2n$) or  
${(N+1)/2}$
elements ($N$ odd, $N=2n+1$),  for example
$\ell^{-\ci}{}_{\!\ci},$
$\ell^{-1}{}_{\!1}$ ... $\ell^{-n}{}_{\!n}$.
For the off-diagonal
$\LLpm $ elements, we can choose as free indices $(C,D)=(c,\ci)$
in  relation (\ref{CLL2}),  and using
$\Lm{\ci}{\ci}\Lm{\bu}{\bu}=\epsi$, we find:
\eq
\Lm{\bu}{c}=
-(C_{\ci \bu})^{-1} C_{ab}\Lm{b}{c}\,\Lm{a}{\ci}\Lm{\bu}{\bu}~.
\label{restrictedL1}
\en
If we choose $(C,D)=(\ci,\ci)$ we obtain
\eq
\Lm{\bu}{\ci}= - (r^{-2} C_{\bu\ci} +C_{\ci\bu})^{-1}
C_{ab}\Lm{b}{\ci}\,\Lm{a}{\ci}\Lm{\bu}{\bu}~. \label{restrictedL2}
\en
Similar results hold for $\Lp{\ci}{d}$ and $\Lp{\ci}{\bu}$.
Iterating this procedure, from
$C_{ab} \Lm{b}{c} \Lm{a}{d} = C_{dc} \epsi$
we find that $\Lm{N}{i}$ (with $i=2,...N-1$) and $\Lm{N}{1}$ are
functionally dependent on $\Lm{i}{1}$ and $\Lm{N}{N}$. Similarly
for $\Lp{1}{i}$ and $\Lp{1}{N}$. The final result is that the  
elements
$~\Lm{a}{J}$ with $\sma{$J<a<J'$}$ and
$~\Lp{a}{J}$ with $\sma{$J'<a<J$}$ --
whose number in both $\pm$ cases is ${1\over 4}N(N+2)$ for $N$ even
and ${1\over 4}(N+1)^2$ for $N$ odd --
and the elements
$\ell^{-\ci}{}_{\!\ci},$
$\ell^{-1}{}_{\!1}$... $\ell^{-n}{}_{\!n}$
generate all $\hat{U}_{q,r}(so(N+2))$.
The total number of  generators is therefore
$(N+2)(N+1)/2$.
\sk

Notice that in this derivation we have not used the $RLL$ relations
(i.e. the quantum analogue of the Lie algebra  relations)
to  further reduce the number of generators. We  therefore  expect  
that,
as in
the classical case,
monomials in the $(N+2)(N+1)/2$ generators can be ordered
(in any arbitrary way).
We begin by proving this for  polynomials  in $\Lp{A}{A}$,
$\Lp{\al}{J}$ with
$\sma{$J'<\al<J$}$,
and
for polynomials in  $\Lm{A}{A}$,
$\Lm{\al}{J}$ with $\sma{$J<\al< J'$}$ .
\sk
{\sl Lemma 4.1 }
Consider the $R\LLpm\LLpm$ commutation  relations
\eq
\R{AB}{EF} \Lpm{F}{D} \Lpm{E}{C} = \Lpm{A}{E} \Lpm{B}{F} \R{EF}{CD}
\label{Rell}~.
\en
{}For $\sma{$C\not= D$}$ they close respectively on the subset of the
$\Lp{\al}{J}$ with
$\sma{$J'<\al\leq J$}$
and on the subset of the $\Lm{\al}{J}$ with $\sma{$J\leq\al<J'$}$.
{}For
$\sma{$C=D$}$ they are
equivalent to the $q^{-1}$-plane commutation relations:
\eq
[P_A\sma{$(J'\!-\!J\!+\!1)$}]^{\al\be}_{~\ga\de}\Lpm{\de}{J}\Lpm{\ga}{ 
J}=0
{}~~,~~~
\label{PAJ'J}
\en where $P_A\sma{$(J'\!-\!J\!+\!1)$}$ is the antisymmetrizer in
dimension
$\sma{$J-J'+1$}$
[compare with (\ref{proiett})].
In particular
\eq
\PA{ab}{cd} \Lm{d}{\ci} \Lm{c}{\ci}=0   \label{PLL}
\en
or equivalently
$[(P_{\!A}){}_{{}_{q^{-1}\!,r^{-1}}}]^{ab}_{~cd}\,\Lm{c}{\ci}
\Lm{d}{\ci}=0$
which coincide, for
$r\rightarrow r^{-1}$ and $q\rightarrow q^{-1}$, with the
$N$-dimensional quantum orthogonal plane relations
(\ref{PRTT13}).

\noi {\sl Proof :} $~~~$The proof is a straightforward calculation   
based
on
(\ref{foraa'}) and on upper or lower triangularity of the $R$ matrix
and of the $\LLpm$ functionals.
\cvd
{\sl Lemma 4.2 } $U_{q,r}(so(N))$ is a Hopf subalgebra of $\U$.

\noi {\sl Proof}:  Choosing $SO_{q,r}(N)$ indices as free indices in
(\ref{Rell})
and using upper or lower triangularity of the $\LLpm$ matrices, and
(\ref{Rnonzero}) or (\ref{Rbig}), we find
that only  $SO_{q,r}(N)$ indices appear in (\ref{Rell}); similarly  
for
relations
(\ref{RLpLm})-(\ref{CLL2}), and for the costructures
(\ref{copLpm})-(\ref{coiLpm}).
\cvd

Now we observe that in virtue of the
$R\LLp\!\LLp$  relations the $\LLp$
elements can be ordered; similarly  we can order the $\LLm$ using the
$R\LLm\!\LLm$ relations.
This statement can be proved by induction using that
$U_{q,r}(so(N))$ is a subalgebra of $\U$, and splitting the $\SO$
index in the usual way [some of the resulting formulas are given in
(\ref{1elle})-(\ref{4elle})].

It is then  straightforward to prove that the elements $\Lp{\al}{J}$  
with
$\sma{$J'<\al\leq J$}$  can
be ordered;
indeed
we can always order the $\Lp{\al}{J}~\Lp{\be}{K}$ with
$\sma{$J'<\al\leq J$}$,
$\sma{$K'<\be\leq K$}$ and $\sma{$J\not= K$}$ since
their commutation relations are a
closed subset
of (\ref{Rell}) [see {\sl Lemma 4.1}].
Then there is no difficulty in ordering substrings
composed by
$\Lp{\al}{J}$ and
$\Lp{\be}{J}$ elements because
(\ref{PAJ'J}) are $q^{-1}$-plane
commutation relations, that allow for any
ordering of the  quantum plane coordinates
 \cite{FRT}.
More in general the  $\Lp{A}{A}$
and $\Lp{\al}{J}$ with
$\sma{$J'<\al<J$}$
can be ordered because of
$
\Lp{A}{A}\Lp{B}{C}$=$(q_{{}_{BA}}/q_{{}_{CA}})\Lp{B}{C}
\Lp{A}{A} \,.$
Similarly we can order the $\Lm{A}{A}$
and $\Lm{\al}{J}$ with
$\sma{$J<\al<J'$}$.
It is now easy to prove the following
\sk
\noi {\sl Theorem 4.1}
A set of elements spanning  $\hat{U}_{q,r}(so(N+2))$ is given by
the ordered monomials
\eq
Mon(\Lp{\al}{J};\sma{$J'<\al<J$})\;(\ell^{-\ci}{}_{\!\ci})^{p_{\ci}}
(\ell^{-1}{}_{\!1})^{p_1}\ldots(\ell^{-n}{}_{\!n})^{p_n}
\;Mon(\Lm{\al}{J};\sma{$J<\al<J'$})
\label{Starrr}
\en
where $p_{\ci}, p_1,...p_n\in {\bf \mbox{\boldmath$N$}}\cup\{0\}$,
$n=N/2$ ($N$ even), $n=(N-1)/2$ ($N$ odd) and
$Mon(\Lp{\al}{J};\sma{$J'<\al<J$})$, $ [Mon(\Lm{\al}{J};
\sma{$J<\al<J'$})]$ is
a monomial in the off-diagonal elements $\Lp{\al}{J}$ with
$\sma{$J'<\al<J$}$
[$\Lm{\al}{J}$ with $\sma{$J<\al<J'$}$] where an ordering has been  
chosen.
\cvd
{\sl Note 4.4 } Conjecture: the above monomials are linearly  
independent
and therefore  form a basis of
 $\hat{U}_{q,r}(so(N+2))\,.$

\sect{Universal enveloping algebra  $U_{q,r}(iso(N))$}
Consider a generic functional $f\in\U$. It is well defined on the
quotient
$\ISO =\SO /H$ if and only if $f(H)=0$.
It is  easy to see that the set $\H$ of all these functionals is a
subalgebra
of $\U$ :
if $f(H)=0$ and $g (H)=0$ then $fg(H)=0$ because
$\Delta  (H) \subseteq
H \otimes \SqrNt + \SqrNt \otimes H.$
Moreover $\H$ is a Hopf subalgebra of $\U$  since $H$ is a Hopf ideal
\cite{Sweedler}.
In agreement with these observations we will find the Hopf algebra
$\IU$ [dually paired to $\ISO$] as a subalgebra of $\U$ vanishing
on the ideal $H$.

Let
\eq
IU\equiv [L^{-A}{}_B, L^{+a}{}_b, L^{+\circ}{}_{\circ},
L^{+\bullet}{}_{\bullet}, \epsi]
\subseteq\U\label{IU}
\en
be the subalgebra of $\U$ generated by
$L^{-A}{}_B, L^{+a}{}_b, L^{+\circ}{}_{\circ},
L^{+\bullet}{}_{\bullet}, \epsi
.$
\sk
{\sl Note 5.1 }:  These are all and only the functionals
annihilating
the generators of $H$:
$\T{a}{\circ}\;,~\T{\bullet}{b}$ and  $\T{\bullet}{\circ}\;$.
The remaining $\U$ generators
$L^{+\circ}{}_b~,~L^{+a}{}_{\bullet}~,~L^{+\circ}{}_{\bullet}$ do not
annihilate
the generators of $H$ and are not included in
(\ref{IU}).
\sk
We now proceed to study  this algebra $IU$. We will show that it is a
Hopf algebra and that $IU\subseteq \H$; we will give an $R$-matrix
formulation,
and prove that $IU$
is a free  $U_{q,r}(so(N))$--module. This is the analogue of
  $ISO_{q,r}(N)$ being a free $SO_{q,r}(N)$-module.
We then show that $IU$ is dually paired with $\ISO$. These results
lead to the conclusion that $IU$ is the
universal enveloping algebra of $ISO_{q,r}(N)$.
\sk
\noi {\sl Theorem 5.1}
$IU$ is a Hopf subalgebra of $\U$.

\noi {\sl Proof :} $~~~$
$IU$ is by definition a subalgebra. The sub-coalgebra property
$\Delta '(IU)
\subseteq IU\otimes IU$
follows immediately from the upper triangularity of $L^{+A}{}_B$:
\eq
\Delta ' (L^{+a}{}_b)=L^{+a}{}_c\otimes L^{+c}{}_b~;~
\Delta '(L^{+\circ}{}_{\circ})=L^{+\circ}{}_{\circ}
\otimes L^{+\circ}{}_{\circ}~;~
\Delta '(L^{+\bullet}{}_{\bullet})
=L^{+\bullet}{}_{\bullet}\otimes L^{+\bullet}{}_{\bullet}
\en
and the compatibility of $\Delta '$ with the product.
We conclude that $IU$ is a Hopf-subalgebra because $\kp(IU)\subseteq  
IU$
as is easily seen using (\ref{coiLpm}) and
antimultiplicativity of $\kp$.
{\cvd}
We may wonder whether the $RLL$ and $CLL$ relations  of $\U$ close in
$IU$.
In this case $IU$ will be given by
all and {\sl only} the polynomials in the functionals
$L^{-A}{}_B, L^{+a}{}_b, L^{+\circ}{}_{\circ},
L^{+\bullet}{}_{\bullet}, \epsi .$
This check is done by writing explicitly all $q$-commutations between
the
generators of $IU$: they do not involve the functionals
$L^{+\circ}{}_b~,~L^{+a}{}_{\bullet}~,~L^{+\circ}{}_{\bullet}$ .
Moreover one can also write them in a compact form using a new
$R$-matrix
${\cal{R}}_{12}\equiv\Lc_2({{t}}_1)$, where $\Lc$ is defined below.
Similarly
the orthogonality
conditions (\ref{CLL1})-(\ref{CLL2})
do not relate elements of $IU$ with elements not belonging to $IU$.
We therefore conclude
\sk
\noi {\sl Theorem 5.2 }
The Hopf algebra $IU$ is generated by the unit $\epsi$ and
the matrix entries:
\eq
L^-=\left(L^{-A}{}_{B_{{}_{}}}\right)
{}~;~~
\Lc=\Mat{L^{+\circ}{}_{\circ}}{0}{0}{0}{L^{+a}{}_b}{0}{0}{0}
{L^{+\bullet}{}_{\bullet}}~;
\en
these functionals satisfy the $q$-commutation relations:
\eq
R_{12} \Lc_2 \Lc_1=\Lc_1 \Lc_2 R_{12} ~\mbox{ or equivalently }~
{\cal{R}}_{12}
\Lc_2 \Lc_1=\Lc_1 \Lc_2 {\cal{R}}_{12} \label{iRLcLc}
\en
\eq
R_{12} \LLm_2 \LLm_1=\LLm_1 \LLm_2 R_{12}~, \label{iRLL}
\en
\eq
{\cal{R}}_{12} \Lc_2 \LLm_1=\LLm_1 \Lc_2 {\cal{R}}_{12}~,
\label{iRLpLm}
\en
where
$$
{\cal{R}}_{12}\equiv\Lc_2({{t}}_1) ~~~
\mbox{ that is }~~~
{\cal{R}}^{ab}_{cd}={R}^{ab}_{cd}~;~~
{\cal{R}}^{AB}_{AB}={R}^{AB}_{AB} ~\mbox{ and otherwise
}~{\cal{R}}^{AB}_{CD}
=0
$$
and the orthogonality conditions :
\eq
C^{AB} \Lc^{C}{}_{B} \Lc^{D}{}_{A} = C^{DC} \epsi ~;~~
C_{AB} \Lc^{B}{}_{C} \Lc^{A}{}_{D} = C_{DC} \epsi ~;\label{iCLcLc}
\en
\eq
C^{AB} \Lm{C}{B} \Lm{D}{A} = C^{DC} \epsi ~;~~
C_{AB} \Lm{B}{C} \Lm{A}{D} = C_{DC} \epsi~, \label{iCLL}
\en
The costructures are the ones given in (\ref{copLpm})-(\ref{coiLpm})
with $\LLp$ replaced by $\Lc$. {\cvd}
{\sl Note 5.2}
We can consider the extension
$\hat{IU}\subset\hat{U}_{q,r}(so(N+2))$ obtained by including
the elements $\ell^{\pm A}{}_A$ (${\ell^{\pm}}^A{}_A=\ln\Lpm{A}{A}$,
see the previous section). Then
$\hat{IU}$
is generated by the symbols $\Lm{A}{B},\, \Lc^A{}_B,\, \ell^{\pm  
A}{}_A$
modulo the relations
(\ref{iRLcLc})-(\ref{iCLL}) and (\ref{1effe1}) [({\ref{epsiaepsi}) in  
the
uniparametric case].
Equivalently,  from (\ref{restrictedL1})-(\ref{restrictedL2}),
we have that $\hat{IU}$ is generated by
$\hat{U}_{q,r}(so(N))$,
the $N$ elements $\Lm{a}{\ci}$ (satisfying  the quantum plane  
relations)
and the dilatation $\ell^{-\ci}{}_{\!\ci}$.
All the relations are then given by  those between the
generators of  $\hat{U}_{q,r}(so(N))$
--listed in  (\ref{RLL})-(\ref{CLL2}), (\ref{1effe1}) with lower case
indices--
and by the following ones
\eqa
& &\Lm{\ci}{\ci}\Lm{a}{\ci}=q^{-1}_{\ci a}\Lm{a}{\ci}\Lm{\ci}{\ci}
\label{1elle}\\
& &P_A^{ab}{}_{\!fe}\Lm{e}{\ci}\Lm{f}{\ci}=0\label{2elle}\\
&
&\Lm{\ci}{\ci}\Lpm{a}{b}={q_{b\ci}\over{q_{d\ci}}}\Lpm{b}{d}\Lm{\ci}{\ 
ci}
\label{3elle}\\
& &\Lm{a}{\ci}\Lpm{b}{d}={r\over {q_{d\ci}}}(R^{\pm})^{ba}{}_{\!ef}
\Lpm{e}{d}\Lm{f}{\ci}\label{4elle}
\ena
where $R^{\pm}$ is defined in (\ref{Rplus}).
The number of generators is $N(N-1)/2+ N + 1$.
\sk
{\sl Note 5.3 }:
When  $q_{a\ci}=r\;\forall a$, then   
$\Lm{\ci}{\ci}=\Lp{\bu}{\bu}\,,\,
\Lm{\bu}{\bu}=\Lp{\ci}{\ci}$ and,  in complete analogy to   
(\ref{(3.24)}),
$IU$ is generated by ${U}_{q,r}(so(N))$,
$\Lm{a}{\ci},\Lm{\ci}{\ci}$ and $\Lm{\bu}{\bu}=(\Lm{\ci}{\ci})^{-1}$.
With abuse of notations we will consider $IU$ generated by these  
elements
also
for arbitrary values of the parameters $q_{a\ci}$; this is  what  
actually
happens in $\hat{IU}$.
\sk
{\sl Note 5.4 }: From the second equation in (\ref{iRLcLc}) applied
to ${{t}}$ we obtain the quantum Yang-Baxter equation for the matrix
${\cal{R}}$.
\sk
{}Following {\sl Note 3.3}, using
(\ref{1elle}), (\ref{2elle})  [quantum plane relations]
and then (\ref{3elle}) and (\ref{4elle}),
a generic element  of $IU$
can be written as $\eta^ia_i$ where $a_i\in U_{q,r}(so(N))$ and  
$\eta^i$
are
the ordered monomials:
$\eta^i=(\Lm{\ci}{\ci})^{i_{\ci}}(\Lm{1}{\ci})^{i_{1}}...\,
(\Lm{N}{\ci})^{i_{N}}$ with $i_{\ci}
\in{\bf \mbox{\boldmath$Z$}}$,
$i_1,... i_N\in
{\bf \mbox{\boldmath$N$}}\cup\{0\}$.
Therefore $IU$ is a  right $U_{q,r}(so(N))$--module
generated by the ordered monomials $\eta^i.$
We now show that as in the classical case the expressions
$\eta^ia_i$
are unique:
$\eta^ia_i=0 \Rightarrow a_i=0 \;\forall\;\sma{$i$}$, i.e.
that $IU$ is a free right $U_{q,r}(so(N))$--module.
To prove this assertion we show that,
at least when  $q_{a\ci}=r\;\forall a$,
$IU$ is a bicovariant bimodule over $U_{q,r}(so(N))$.
Since any bicovariant bimodule is free\footnote{The results of  
\cite{Wor}
apply to a general Hopf algebra with invertible antipode. This can be
shown by checking that
all the formulae used to derive the results of \cite{Wor} --they are
collected in the appendix
of \cite{Wor}-- hold also in the general case of a Hopf algebra with
invertible antipode.}
 \cite{Wor} we then deduce that,  as a right module, $IU$
is freely generated by the $\eta^i$.
\sk
\noi {\sl Theorem 5.3 } Consider $IU$  (with the parameter  
restriction
$q_{a\ci}\!=r\;\forall a$) as the right $U_{q,r}(so(N))$--module
$\Ga=\{\eta^ia_i\}$ [$a_i\in U_{q,r}(so(N))$] generated by the  
ordered
monomials
$\eta^i=(\Lm{\ci}{\ci})^{i_{\ci}}(\Lm{1}{\ci})^{i_{1}}...\,
(\Lm{N}{\ci})^{i_{N}}$
with $i_{\ci}\in
{\bf \mbox{\boldmath$Z$}}$,
$i_1,... i_N\in {\bf \mbox{\boldmath$N$}}\cup\{0\}$ .
\begin{description}
\item[a)]{$\Ga$ is a bimodule with the left module structure  
trivially
inherited from the algebra $IU$}.
\item[b)]{$\Ga$ is a right covariant bimodule with right coaction
$\de_R\,:\;\Ga\rightarrow \Ga \otimes U_{q,r}(so(N))$ defined by
\eq
\de_R(\eta^i)\equiv\eta^i\otimes \epsi ~~,~~~
\de_R(a\eta^i b)\equiv\D '(a)\de_R(\eta^i)\D '(b)~.
\label{dercoa}
\en
}
\item[c)]{$\Ga$ is a left covariant bimodule with left coaction
$\de_L\;:~ \Ga\rightarrow U_{q,r}(so(N))\otimes \Ga$ defined by
\eqa
&\de_L(\Lm{\ci}{\ci})\equiv \epsi\otimes\Lm{\ci}{\ci}~;~~
\de_L(\Lm{a}{\ci})\equiv\Lm{a}{b}\otimes\Lm{b}{\ci}&\\
&\de_L(a\Lm{\al}{\ci}\Lm{\be}{\ci}... \Lm{\ga}{\ci}b)\equiv
\D '(a)\de_L(\Lm{\al}{\ci})\de_L(\Lm{\be}{\ci})...  
\de_L(\Lm{\ga}{\ci})\D
'(b)
\ena
where $\sma{$\al=(\ci,a)$}$,  $\sma{$\be=(\ci,b)$}$,
$\sma{$\ga=(\ci,c)$}$.
}
\item[d)]{$\Ga$ is a bicovariant bimodule
\eq
(id\otimes\de_R)\de_L=(\de_L\otimes id)\de_R~.
\label{Bicco}
\en
}
\item[e)]{$\Ga$ is freely generated by the right invariant elements
$\eta^i~.$}
\end{description}
\noi {\sl Proof :} $~~~$

\noi{\bf a)}{} is immediate since, from {\sl Note 5.3} and
{\sl Lemma 4.2}, $U_{q,r}((so(N))$
is a subalgebra of $IU$.

\noi{\bf b)}{} Consider the linear map
$\de_r : IU\rightarrow IU\otimes IU$
defined by
\eq
\de_r(\Lm{\al}{\ci})=\Lm{\al}{\ci}\otimes \epsi~;~~\de_r(a)=\D '(a)~
\forall\,a\in U_{q,r}(so(N))
\label{derIU}
\en
and extended multiplicatively on all $IU$.
This map is obviously well defined on
$U_{q,r}(so(N))$ because it coincides with the coproduct on
$U_{q,r}(so(N))$
[$U_{q,r}(so(N))$ is a Hopf subalgebra of $IU$]; it is also well  
defined
on
all $IU$ since it is multiplicative and compatible with
(\ref{1elle})-(\ref{4elle}). We check for example
(\ref{4elle}) with
$q_{a\ci}=r\;\forall a$:
$$
\de_r(\Lm{a}{\ci}\Lpm{b}{d})=
\Lm{a}{\ci}\Lpm{b}{c}\otimes\Lpm{c}{d}=
(R^{\pm})^{ba}{}_{\!ef}\Lm{e}{c}\Lm{f}{\ci}\otimes
\Lpm{c}{d}=\de_r\!\left( (R^{\pm})^{ba}{}_{\!ef}
\Lpm{e}{d}\Lm{f}{\ci}\right) .
$$
This shows that $\de_R\;:~\Ga\rightarrow \Ga\otimes U_{q,r}(so(N))$
is well defined
since $\Ga$ is  $IU$ seen as a $U_{q,r}(so(N))$-bimodule and the  
actions
of
$\de_r$ and $\de_R$ on $\Ga$ coincide.

It is now  immediate to show
that $\Ga$ is a right covariant bimodule, i.e. that
\eq
\forall \;\eta^ia_i\in \Ga ;~~~(\de_R\otimes id)\de_R(\eta^ia_i)=
(id\otimes \D ')\de_R(\eta^ia_i)~;~~(id\otimes\epsi  
')\de_R(\eta^ia_i)=
\eta^ia_i~.
\en
{\bf c)}{} We proceed as in the previous case, defining the
linear map  $\de_l:IU\rightarrow IU\otimes IU$,
\eq
\de_l(\Lm{a}{\ci})=\Lm{a}{b}\otimes \Lm{b}{\ci}~;~~
\de_l(\Lm{\ci}{\ci})=\Lm{\ci}{\ci}\otimes  
\Lm{\ci}{\ci}~;~~\de_l(a)=\D
'(a)~
\forall\,a\in U_{q,r}(so(N))
\label{delIU}
\en
which is  extended multiplicatively on all $IU$.
As was the case for $\de_r$, it is well defined on
$U_{q,r}(so(N))$ and it is also well defined on all $IU$ because it
is multiplicative and compatible with
(\ref{1elle})-(\ref{4elle}). For example, the compatibility of   
$\de_l$
with relation (\ref{2elle})
 holds because
$P_A^{ab}{}_{\!ef}\Lm{f}{d}\Lm{e}{c}=\Lm{b}{f}\Lm{a}{e}P_A^{ef}{}_{\!c 
d}$
[a consequence of $(\hat{R})^{\pm 1}\LLpm_2\LLpm_1=
\LLpm_2\LLpm_1(\hat{R})^{\pm 1}$ and the fact that $P_A$ is a   
polynomial
in
$\hat{R}$ and $\hat{R}^{-1}$, see (\ref{proiett})]. This is in  
complete
analogy with the compatibility  of the
 left coaction $\de(x^a)=\T{a}{b}\otimes x^b$ with the $q$-plane
commutation relations.

To prove that $\Ga$ is a left covariant bimodule, notice that
\eq
(\epsi\otimes id)\de_l(\Lm{a}{\ci})=\Lm{a}{\ci}\;,~
(\D '\otimes
id)\de_l(\Lm{a}{\ci})=\Lm{a}{d}\otimes\Lm{d}{b}\otimes\Lm{b}{\ci}
=(id\otimes\de_l)\de_l(\Lm{a}{\ci})~,
\en
and similarly for $\Lm{\ci}{\ci}$. Now since $\de_l(a)=\D '{(a)}$ if
$a\in U_{r}(so(N))$, and since $\de_l$ is
multiplicative, we have on all $IU$
\eq
(\epsi\otimes id)\de_l=id~;~~(\D '\otimes id)\de_l=(id\otimes  
\de_l)\de_l
\en
{\bf d)}{} The bicovariance condition (\ref{Bicco}) follows directly
from:
\eqa
&(id\otimes  
\de_r)\de_l(\Lm{a}{\ci})=\Lm{a}{b}\otimes\Lm{b}{\ci}\otimes
\epsi=
(\de_l\otimes id)\de_r(\Lm{a}{\ci})&\\
&(id\otimes \de_r)\de_l(\Lm{\ci}{\ci})= \epsi\otimes\Lm{\ci}{\ci}
\otimes \epsi=(\de_l\otimes id)\de_r(\Lm{\ci}{\ci})&
\ena
{\bf e)}{} We now recall that a bicovariant bimodule is always
freely generated by a basis of $\Ga_{inv}$, the space of
right invariant elements of
$\Ga$ \cite{Wor}.
We also know that the $\eta^i$
are right invariant. Now, since they generate $\Ga$, they linearly  
span
$\Ga_{inv}$, and since they are linearly independent
, they form a basis of $\Ga_{inv}$.
We conclude that $\Ga$ is freely generated by the $\eta^i$:
$\eta^ia_i=0\Rightarrow a_i=0~\forall \sma{$i$}$.
\cvd
\nopagebreak
It is now  easy to prove that the
$\eta^i$ freely generate $IU$ also without the restriction
$q_{a\ci}=r\;\forall a$. [Hint: recall the definition of $\LLm$
as:
${\Lm{A}{B}(c)=\cal{R^{(F)}}}^{-1}(\T{A}{B}
\otimes {c})~\forall\,c\in\SO$, and use
${\cal{F}}\in \hat{U}^0\otimes\hat{U}^0$
to show that $\Lm{A}{B}$
differs from the uniparametric $\Lm{A}{B}$
(obtained with $\cal{R}$ instead of ${\cal{R^{(F)}}}$)
by a factor
belonging to $\hat{U}^0$ and invertible.]
\sk
\vskip .2cm
\noi {\large{{\bf{Duality}} $\IU\leftrightarrow \ISO$}}
\sk
We now show that $IU$ is dually paired to $\SO$. This is the
fundamental
step allowing to interpret $IU$ as the algebra of regular functionals
on $\ISO$.
\sk
\noi {\sl Theorem 5.4} $~IU$ annihilates $H$, that is $IU\subseteq
\H$.

\noi {\sl Proof :} $~~~$
Let ${\cal{L}}$ and $\Tc$ be generic generators of $IU$ and  $H$
respectively.
As discussed in {\sl Note 5.1}, ${\cal{L}}(\Tc)=0$. A generic element
of the ideal is given by $a\Tc b$ where sum of polynomials is  
understood;
we have (using Sweedler's notation for the coproduct):
$
{\cal{L}}(a\Tc
b)={\cal{L}}_{(1)}(a){\cal{L}}_{(2)}(\Tc){\cal{L}}_{(3)}(b)=0
$
because ${\cal{L}}_{(2)}(\Tc)=0$. Indeed ${\cal{L}}_{(2)}$ is
still a generator of $IU$ since $IU$ is a sub-coalgebra of $\U.$
Thus ${\cal{L}}(H)=0$. Recalling that a product of functionals
annihilating
$H$ still annihilates the co-ideal $H$, we also have $IU(H)=0$.
{\cvd}
In virtue of {\sl Theorem 5.4} the following bracket is well defined:
\eqa
\!\!\!\!\!\!\!\!\!\!\!\!\!\!\!\!\!\!\mbox{{\sl Definition }$~~~~~$}
&\!\!\!\!\!\!\!\!\!\le ~ ,~ \re \; :~ &IU\otimes \ISO \longrightarrow
\mbox{\boldmath $C$} \nonumber\\
&               &\le a', P(a)\re\equiv a'(a)
 ~~~\forall \/a'\in IU \,,~\forall \/a\in \SO\label{duality}
\ena
where $P~:~\SO\rightarrow \SO/H\equiv \ISO$ is the canonical
projection, which is surjective.
The bracket is well
defined because two generic counterimages of $P(a)$ differ
by an addend belonging to $H$.

Note that when we use the bracket $\le ~,~\re$, $a'$ is seen as an
element
of $IU$ , while in the expression $a'(a)$, $a'$ is seen
as an element of $\U$ (vanishing on $H$).
\sk
\noi { \sl Theorem 5.5 $~~~$} The bracket (\ref{duality})
defines a pairing between $IU$ and $\ISO$ :
$\forall\/ a',b'\in IU~,~\forall\/P(a),P(b)\in \ISO$
\eqa
& &\le a'b' , P(a)\re = \le a'\otimes b',\Delta
(P(a))\re\label{uuno}\\
& &\le a',P(a)P(b)\re=\le\Delta '(a'),P(a)\otimes
P(b)\re\label{udue}\\
& &\le\kp(a'),P(a)\re=\le a',\kappa(P(a))\re\label{utre}\\
& &\le I,P(a)\re=\epsi (P(a))~~;~~~\le a',P(I)\re=\epsi
'(a')\label{uquattro}
\ena

\noi{\sl Proof : } The proof is easy since $IU$ is a Hopf subalgebra
of $\U$ and $P$ is compatible with the structures and costructures
of $\SO$ and $\ISO$. Indeed we have
\[
\le a',P(a)P(b)\re=\le a',P(ab)\re=a'(ab)=\Delta '(a')(a\otimes b)
=\le\Delta ' (a'), P(a)\otimes P(b)\re
\]
\[
\le a'b',P(a)\re=a'b'(a)=(a'\otimes b')\Dtwo (a)=
\le a'\otimes b',(P\otimes P)\Dtwo(a)\re=\le a'\otimes b',\Delta
(P(a))
\re
\]
\[
\le \kp(a'),P(a)\re=\kp(a')(a)=a'(\kappa_{N+2}(a))
=\le a',P(\kappa_{N+2}(a))\re=\le a' ,\kappa(P(a))\re
\]
\sk
\cvd
We now recall that $IU$ and $\ISO$, besides being dually paired,
are free right modules respectively on $U_{q,r}(so(N))$
and  on $SO_{q,r}(N)$.
They are  freely
generated by the two isomorphic sets of the ordered monomials in
 $\Lm{\ci}{\ci},\;\Lm{a}{\ci}$
and $u,\;x^a$ respectively [cf. the commutations (\ref{1elle}),
(\ref{2elle}) and (\ref{PRTT22}),
(\ref{PRTT13})].  We can then call
$IU$  the universal enveloping algebra of  $\ISO$
\eq
\IU\equiv IU~
\en
in the same way  $U_r(so(N))$ is
the universal enveloping algebra of $SO_r(N)$ \cite{FRT}.
\sk
\noi{\sl Note 5.5 : } Given a $*$-structure on $\ISO$, the duality
$\ISO\leftrightarrow \IU$
induces a $*$-structure on $\IU$. If in particular the
$*$-conjugation
on $\ISO$ is
found by
projecting a $*$-conjugation on $\SO$, then the induced $*$ on $\IU$
is simply the
restriction to $\IU$ of
the $*$ on $\U$. This is the case for the $*$-structures that lead
to the real
forms
$ISO_{q,r}(N, \mbox{\boldmath $R$})$ and $ISO_{q,r}(n+1, n-1)$ and in
particular to
the quantum Poincar\'{e} group \cite{Firenze1,Cas2,inson}.

\sect{Projected differential calculus}
\sk
In the previous sections
we have found the inhomogeneous
quantum group  $\ISO$ by means of a projection from
$\SO$.
Dually, its universal
enveloping algebra is a given  Hopf subalgebra of
$\U$.
Using the same techniques  differential calculi on
$\ISO$ can be found.
\sk
\noi {\large {\bf Projecting Woronowicz ideal}}
\sk
Following Woronowicz \cite{Wor},  we recall that a bicovariant
differential calculus
over a generic Hopf algebra $A$ is determined by a right ideal $\Re$
of $A$. This ideal  has to be included in ker$\epsi$ (i.e.
its elements have vanishing counit) and must be ad-invariant that is,
$a\!d_A(r) \in \Re\otimes A ~~\forall r\in \Re $
 where $a\!d_A(r)$ is defined by
$a\!d_A(a)\equiv a_2\otimes\kappa_A(a_1)a_3~~~\forall\/a\in A~~ $;
the
index {{\scriptsize $ A$}} denotes the costructures in $A$
and we have used Sweedler's notation for the coproduct.
For any such  $\Re$ one can construct a bicovariant differential
calculus.
In the following we show that , given a quotient Hopf algebra
$A/H$ (with  canonical projection $P:A\rightarrow A/H \equiv P(A)$),
 $P(R)$ is a right
ad-invariant ideal in  $P(A)$; therefore it defines
a bicovariant differential calculus at the projected level.
Moreover the space of tangent vectors on $P(A)$ is easily found as a
subspace
of
the tangent vectors on $A$. The explicit construction of the exterior
differential
$d$, and of the bicovariant bimodule $\Gamma$ of one-forms
is then
straightforward.
\sk
\noi {\sl Theorem 6.1}$~~~$
If $\Re\in $ ker$\epsi$ is a right ad-invariant ideal of $A$
then $P(\Re)$ is included in ker$\epsi$ and is a right ad-invariant
ideal of $P(A).$

\noi {\sl Proof : }
The only nontrivial part is ad-invariance.
{} From $ a\!d_A(r)=r_2\otimes\kappa_A(r_1)r_3\in \Re\otimes
A~~\forall r\in
\Re~ $,
applying $P\otimes P$ we obtain
$P(r_2)\otimes\P(\kappa_A(r_1))P(r_3)\in
P(\Re)\otimes P(A)$ $\forall P(r)\in P(\Re).$
Now
\eq
P(r_2)\otimes\P(\kappa_A(r_1))P(r_3)=P(r_2)\otimes
\kappa(P(r_1))P(r_3)=
P(r)_2\otimes\kappa(P(r)_1)P(r)_3 \equiv a\!d(P(r))\label{Pad}
\en
where we have used compatibility of the projection with the
costructures of $A$ and $P(A)$; $\kappa$ denotes the antipode in
$P(A)$ and,
after the second equality, the coproduct of $P(A)$ is understood.
Relation (\ref{Pad}) gives the ad-invariance of $P(\Re)$:
$
\forall P(r)\in P(\Re)$ $a\!d(P(r))\in P(\Re)\otimes P(A).
$
{\cvd}
The space of tangent vectors on a quantum group $P(A)$ is given by
\cite{Wor}:
\eq
T\equiv\{\bar{\chi} : P(A)\rightarrow \mbox{\boldmath $C$} ~~ |~
\bar{\chi}
\mbox{ linear functionals, }
\bar{\chi}(I)=0 \mbox{ and } \bar{\chi}(P(\Re))=0\}~.\label{tangenti}
\en
{\sl Remark:} the vector space $T$ defined in
(\ref{tangenti})
is given
by all and only those functionals $\bar{\chi}$
corresponding to elements $\chi$ of the tangent space  $T_{A}$ on $A$
that annihilate the Hopf ideal $H$.
Indeed
if $\chi$ annihilates $H$,
then
$\bar{\chi}$ defined by
$
\bar{\chi} : {A/H} \longrightarrow \mbox{\boldmath $C$}$ with
$\bar{\chi}(P(a)) \equiv \chi(a),~
\forall P(a)\in P(A)$
is a well defined functional on $P(A)$
[see (\ref{duality})].
{}From $\chi(\Re)=0$ we obtain
$\bar{\chi}(P(\Re))= 0 $
i.e. $\bar{\chi}\in T$. Viceversa a functional
$\bar{\chi}\in T$ is trivially extended to a functional $\chi\in  
T_A$.
\sk
Recall \cite{Wor,SWZ} that the deformed Lie bracket
is given by  $[\chi_i,\chi_j](a)=(\chi_i\otimes\chi_j)a\!d_A(a)$
where $\chi_i , \chi_j $ are functionals on $A$.
For the ``projected'' $q$-Lie algebra we have:
\sk
\noi {\sl Theorem 6.2$~~~$} The $q$-Lie algebra on $P(A)$ is a closed
subset
of the $q$-Lie algebra on $A$.

\noi {\sl Proof : }
Let $\chi_i(H)=\chi_j(H) =0$. We have, using (\ref{Pad}) in the
second equality
$$
[\bar{\chi}_i,\bar{\chi}_j](P(a))=
(\bar{\chi}_i\otimes \bar{\chi}_j)a\!d(P(a))=
\bar{\chi}_i\otimes \bar{\chi}_j(P\otimes P)a\!d_A(a)=
(\chi_i\otimes\chi_j)a\!d_A(a)=[\chi_i,\chi_j](a)
$$
in particular
$[\bar{\chi}_i,\bar{\chi}_j](P(\Re))=[{\chi}_i,{\chi}_j](\Re)=0$
and this proves the theorem.
\cvd
\noindent In virtue of {\sl Theorem 6.2}  the following corollary is
easily
proved.
\sk
{\sl Corollary }  Consider the structure constants
$\mbox{\boldmath $C_{ij}{}^{k}$}\!$
defined by $[\chi_i ,\chi_j] =
\mbox{\boldmath $C_{ij}{}^{k}\!\!$}~ \chi_k$,
where  $\{\chi_i\}$ will henceforth denote  a basis of $T_A$  
containing
the
maximum number of tangent vectors vanishing
on $H$. The subset of
the structure constants corresponding to the functionals  $\chi_i$
that
annihilate $H$
is the set of all the structure constants of $P(A)$.
\cvd
The exterior differential related to this projected calculus is given  
by:
\eq
\forall \/a\in P(A)~~~~~~da=(\bar{\chi}_i*a)\bar{\omega}^i
\label{dachia}
\en
where $\bar{\chi}_i*a\equiv (id\otimes \chi_i)\Delta a$, and
$\bar{\omega}^i$
are the one-forms dual to the tangent vectors $\bar{\chi}_i$
\cite{Wor,PaoloPeter}; they freely generate the left module of
one-forms
$\Gamma=\{a_i\bar{\omega}^i~,~~a_i\in P(A)\}$. The right module
structure is
given by
the $\bar{f}^i{}_j$ functionals, obtained applying the coproduct
$\D'$ to
the
$\bar{\chi}_i$:
\eq
\D'\bar{\chi}_i=\bar{\chi}_j\otimes \bar{f}^j{}_i +
\epsi\otimes\bar{\chi}_i
{}~~~\Rightarrow~~~\bar{\omega}^ia=(\bar{f}^i{}_j*a)\bar{\omega}^j~.
\label{rightmodule}
\en
The space $\Ga$ of one-forms on $P(A)$  can be studied by projecting  
the
one-forms on
$A$ into one-forms on $P(A)$. For this we introduce the projection  
$P$
acting on
$\Ga_A$ (the space of one-forms on $A$) as follows:
\eqa
\!\!\!\!\!\!\!\!\!\!\!\!\!\!\!\!\!\!\!\!\!\!\!\!\!\!\!\!\!\!\!\!\!\!\! 
\!\!
\!\!\!\!
\mbox{{\sl Definition $~~~~~~~~~~~~~~~~~~~~~~$}}
& P~:~&\Gamma_{\!A}\rightarrow\Gamma\\
& & a_i\omega^i\mapsto P(a_i)\bar{\omega}^i
\ena
where $\bar{\omega}^i =0$ if ${\chi}_i(H)\not= 0$.  We now show that  
$P$
is a
bicovariant bimodule epimorphism and that it is compatible with the
differential
calculi.
Trivially $P$ is a left module epimorphism because $\Gamma_{\!A}$ and
$\Gamma$ are
free left modules generated respectively by the one-forms
$\{\omega^i\}$ and
$\{\bar{\omega}^i\}$. It is also easy to see [use
(\ref{rightmodule})] that
$
\forall\/\rho\in \Gamma_{\!A}$ , $\forall\/a\in A~~$ $P(\rho
a)=P(\rho)P(a)\,,
$
which shows that $P$ is a bimodule epimorphism.

To prove  that $P$ is compatible with the exterior
differentials
$d_{\!A}$ on $A$ and $d$ on $P(A)$, consider the generic one-form
$a\,d_{\!A}b=a(\chi_i*b)\omega^i$
[see (\ref{dachia})]; we have
$P(a \,d_{\!A}b)$=$P(a)P(\chi_i*b)\bar{\omega}^i$=$
P(a)[\bar{\chi}_i*P(b)]\bar{\omega}^i$
=$P(a)dP(b)\,.$

Since the exterior differential $d$ induces the comodule structure on
$\Gamma$ by the definitions:
\eq
\forall\/ a,b\in P(A)~~~~~~~~~~
\begin{array}{l}
\Delta_{\!L}(a d b)\equiv \Delta(a)(id\otimes d)\Delta(b)~,\\
\Delta_{\!R}(a d b)\equiv \Delta(a)(d\otimes id)\Delta(b)~,
\end{array}
\en
Finally  $P$ is a comodule homomorphism:
$\Delta_{\!L}(P(\rho))=(P\otimes P)\Delta_{\!L A}(\rho)\,$;
$\Delta_{\!R}(P(\rho))=(P\otimes P)\Delta_{\!R A}(\rho)$,
$\forall\/\rho\in\Gamma_A$
where $\Delta_{\!L A}~(\Delta_{\!R A})$  is the left (right) coaction
of $A$.

{} From $\Delta_{\!L A}\omega^i=I\otimes\omega^i$ and
$\Delta_{\!R A}\omega^i=\omega^j\otimes M_j{}^i$, where $M_j{}^i$
defines the
adjoint representation on $A$, we obtain an explicit expression for
$\D_L$ and $\D_R$:
\eq
\Delta_{\!L}\/\bar{\omega}^i=I\otimes\bar{\omega}^i~~;~~~
\Delta_{\!R}\/\bar{\omega}^i=\bar{\omega}^j\otimes P(M_j{}^i)~.
\label{adjcoaction}
\en
\sk
\noi {\large{\bf Application: $ISO_{q,r}(N)$ differential calculi}}
\sk
We now apply the above discussion to the quantum groups $A=\SO$ and
$P(A)=\ISO$.
The adjoint representation for
$\SO$ is given by
\eq
\MM{A}{BC}{D}\equiv\T{A}{C}\kappa_{N+2}(\T{D}{B})~,\label{last?}
\en
and the $\chi$
functionals explicitly read
\eq
\cchi{A}{B} = \rinv [\ff{C}{CA}{B}-\de^A_B \epsi]~~
\mbox{ where }~~
\ff{A_1}{A_2B_1}{B_2} \equiv \kp (\Lp{B_1}{A_1}) \Lm{A_2}{B_2},
\label{defff}
\en
see  \cite{Jurco} and references therein (see also
\cite{altroarticolo}).
Decomposing the indices we find:
\eqa
& & \cchi{a}{b}=\rinv [\ff{c}{c a}{b}-\de^a_b \epsi]~~~~~~~+\rinv
\ff{\bu}{\bu
a}{b}\label{chiab} \\
& & \cchi{a}{\ci}=\rinv \ff{c}{c a}{\ci} ~~~~~~~~~ ~~~~~~~ +\rinv
\ff{\bu}{\bu
a}{\ci} \\
& &\cchi{\ci}{b}=~~~~~~~~~~~~~~~~~~~~~~~~~~~~~~~~~~ + \rinv
 [\ff{c}{c \ci}{b}+\ff{\bu}{\bu \ci}{b}]  \\
& & \cchi{a}{\bu}=~~~~~~~~~~~~~~~~~~~~~~~~~~~~~~~~~~+ \rinv
\ff{\bu}{\bu a}{\bu}\\
& & \cchi{\bu}{b}=\rinv \ff{\bu}{\bu\bu}{b} ~~~~~~~~~~~~~\\
& &\cchi{\ci}{\ci}=\rinv [\ff{\ci}{\ci\ci}{\ci}-\epsi]~~~~~~~~~~ +
\rinv
[\ff{c}{c\ci}{\ci}
+\ff{\bu}{\bu \ci}{\ci} ] \\
& &\cchi{\ci}{\bu}=~~~~~~~~~~~~~~~~~~~~~~~~~~~~~~~~~~ + \rinv
\ff{\bu}{\bu \ci}{\bu}\\
& & \cchi{\bu}{\ci}=\rinv \ff{\bu}{\bu \bu}{\ci}~~~~~~~~~~~~~~~\\
& & \cchi{\bu}{\bu} = \rinv [\ff{\bu}{\bu\bu}{\bu} - \epsi]~~~~~\\
& &~~~~~~~~\underbrace{~~~~~~~~~~~~~~~~~~~~~~~~~~~}_
{\hbox{\small{terms annihilating $H$}}} \nonumber
\ena
\noi where using {\sl Theorem 5.4} and {\sl Note 5.1} we have
indicated the
terms that do and do not
annihilate the Hopf ideal
$H$. We see that only three of these functionals, namely
$\cchi{\bu}{b}$,  $\cchi{\bu}{\ci} $ and $\cchi{\bu}{\bu}$,
do vanish on
$H$.
The resulting bicovariant differential
calculus contains dilatations and translations,
but does not contain the tangent vectors of $SO_{q,r}(N)$, i.e.
the functionals $\cchi{a}{b}$.
The differential related to this calculus is given by
\eq
\forall a\in \ISO ~~~~~~da=(\cchi{\bu}{b}*a)\om_{\bu}{}^b
+(\cchi{\bu}{\bu}*a)\om_{\bu}{}^{\bu} +
(\cchi{\bu}{\ci}*a)\om_{\bu}{}^{\ci}\label{dMink}
\en
where $\om_{\bu}{}^{b},~\om_{\bu}{}^{\bu}$ and $\om_{\bu}{}^{\ci}$
are the
one-forms dual
to the tangent vectors $\cchi{\bu}{b}$,  $\cchi{\bu}{\ci} $ and
$\cchi{\bu}{\bu}$
\cite{Wor,PaoloPeter}
(with abuse of notation, we omit the bar over the elements of the
projected
calculus).
The $q$-Lie algebra is explicitly given by\footnote{We thank A.
Scarfone for
the derivation
of (\ref{quattro}).}
\eqa
&&\cchi{\bu}{\ci}\cchi{\bu}{b}-(q_{\bu b})^{-2}
\cchi{\bu}{b}\cchi{\bu}{\ci}=0\\
&&\cchi{\bu}{c}\cchi{\bu}{\bu}-r^{-2}\cchi{\bu}{\bu}
\cchi{\bu}{c}=-r^{-1}\cchi{\bu}{c}\label{qlieMin}\\
&&\cchi{\bu}{\ci}\cchi{\bu}{\bu}-r^{-4}\cchi{\bu}{\bu}
\cchi{\bu}{\ci}={-(1+r^2)\over {r^{3}}}\cchi{\bu}{\ci}\\
&& q_{\bu a} \PA{ab}{cd}\cchi{\bu}{b} \cchi{\bu}{a}=0
\label{quattro}
\ena
A combination of  the above relations yields:
\eq
\cchi{\bu}{\ci} + \lambda \cchi{\bu}{\ci} \cchi{\bu}{\bu} = \lambda
{-r^{{N\over 2}}\over {r^2 + r^N}}
{1\over q_{d\bu}} \cchi{\bu}{b} C^{db} \cchi{\bu}{d}  \label{chibuci}
\en
Notice the similar structure of eq.s (\ref{PRTT44}),  
(\ref{restrictedL2})
and
(\ref{chibuci}).
\sk
\noi The bicovariant bimodule of one-forms is characterized by the
functionals
\eq
\ff{\bu}{\ci \bu}{\ci}~,~~\ff{\bu}{a \bu}{\ci}~,~~\ff{\bu}{\bu
\bu}{\ci}~,~~
\ff{\bu}{a \bu}{b}~,~~\ff{\bu}{\bu \bu}{b}~,~~\ff{\bu}{\bu \bu}{\bu}~
\en
that appear in the comultiplication of
$\cchi{\bu}{b}$,  $\cchi{\bu}{\ci} $ and $\cchi{\bu}{\bu}$
[use upper (lower) triangularity of $L^+ \;(L^-)$],
and  by the elements
\eq
P({\MM{\bu}{B\bu}{D}}) =
P(\T{\bu}{\bu}\/\kappa_{N+2}(\T{D}{B}))=vP(\kappa_{N+2}(\T{D}{B}))
\label{adMink}
\en
that, according to (\ref{last?}) and (\ref{adjcoaction}),
characterize the right coaction of $\ISO$ on
$\om_{\bu}{}^{b},~\om_{\bu}{}^{\bu}$ and $\om_{\bu}{}^{\ci}$.
They explicitly read
\eq
\begin{array}{lll}
P(\MM{\bu}{\ci\bu}{\ci})=v^2     &P(\MM{\bu}{\ci\bu}{d})=0
&P(\MM{\bu}{\ci\bu}{\bu})=0     \\
P(\MM{\bu}{b \bu}{\ci})=vr^{-{N\over 2}}x^eC_{eb}      &P(\MM{\bu}{b
\bu}{d})=v\kappa(\T{d}{b})    &P(\MM{\bu}{b \bu}{\bu})=0     \\
P(\MM{\bu}{\bu\bu}{\ci})=-{1\over{r^N(r^{N\over 2}+r^{-{N\over
2}+2})}}x^eC_{ef}x^f
&P(\MM{\bu}{\bu\bu}{d})=v\kappa(x^d)     &P(\MM{\bu}{\bu\bu}{\bu})=I
{}~~~
\label{cccinque}
\end{array}
\en
Notice that only the couples of indices $\sma{$(\bu \ci),~(\bu b)$}$
and
$\sma{$(\bu\bu)$}$ appear
in (\ref{dMink})-(\ref{cccinque}): these are therefore the only  
indices
involved
in the projected differential calculus on $ISO_{q,r}(N)$.
\sk
The functionals $\cchi{a}{b}$ cannot be good tangent vectors
on $\ISO$ because of the functionals
$\ff{\bu}{\bu a}{b}$ appearing in (\ref{chiab}): these  do not
annihilate $H$.
We see however that
$
\limrone{1\over {r-r^{-1}}}\ff{\bu}{\bu a}{b}(a)=0 ~~~~\forall a\in
\SO\;;
$
for this reason we consider in the following
the particular multiparametric deformations called ``minimal
deformations" (twistings), corresponding to $r=1$.

As shown in \cite{altroarticolo} in the
$r \rightarrow 1$ limit  the  $\chi$
functionals are given by:
\eqa
& &\cchi{A}{A}= \limrone \linv~ [\ff{A}{AA}{A}-\epsi]~~~;~~~~~~~~
\cchi{A}{A'}=0\nonumber\\
& &\cchi{A}{B}=\limrone \linv  
~\ff{A}{AA}{B},~~~~~\sma{$A>B$}~~~;~~~~~~~~
\cchi{A}{B}=\limrone \linv ~ \ff{B}{BA}{B},~~~~~\sma{$A<B$}
\nonumber\ena
where $\lambda\equiv r-r^{-1}$, and close on the $q$-Lie algebra
\eqa
& &\cchi{C_{1}}{C_2}  \cchi{B_{1}}{B_2} - q_{B_1C_2} q_{C_1B_1}
q_{B_2C_1} q_{C_2B_2}
{}~\cchi{B_{1}}{B_2} \cchi{C_{1}}{C_2}=~~~~~~~~~\nonumber\\
& &~~~~-q_{B_1C_2} q_{C_2B_2} q_{B_2B_1} \de^{C_1}_{B_2}
{}~\cchi{B_{1}}
{C_2}+
q_{C_1B_1} q_{B_2B_1} C_{B_2C_2} ~\cchi{B_{1}}{C_1'}+\nonumber\\
& &~~~~+q_{C_2B_2} q_{B_1C_2} C^{C_1B_1}~ \cchi{B_2'} {C_2}-
q_{B_2C_1} \de^{B_1}_{C_2}~\cchi{B_2'}{C_1'}~. \label{Lierone}
\ena
 Not all of these functionals are linearly
independent because:
\eq
\cchi{B'}{A'}=-q_{AB}\,\cchi{A}{B}~.\label{dependence}
\en
{}From (\ref{dependence}) we see that a basis of
tangent vectors on $SO_{q,r=1}(N+2)$ is given by
\eq
\{\cchi{A}{B}~,~\mbox{with }
\sma{$A+B > {N+1}~,~~A,B: 0=\ci,1,2,...,N,N+1=\bu$}\}~.
\en
They define a bicovariant differential calculus on $SO_{q,r=1}(N+2)$.
The projected
bicovariant calculus on $ISO_{q,r=1}(N)$
is therefore
characterized by the basis
of tangent vectors
\eq
\cchi{a}{b}=\limrone \linv~ [\ff{c}{c a}{b}-\de^a_b \epsi]~,~~~
\mbox{with  \sma{$a+b> N+1~$}};\label{ISOtang1}
\en
\eq
\cchi{\bu}{b}=\limrone \linv \ff{\bu}{\bu\bu}{b}~~;~~~
\cchi{\bu}{\bu} = \limrone \linv~ [\ff{\bu}{\bu\bu}{\bu} - \epsi]~,
\label{ISOtang2}
\en
indeed {\sl Theorem 5.4} assures that these functionals
annihilate $H$, while from {\sl Note 5.1} it is not difficult to see
that
the remaining functionals $\cchi{a}{\bu} = \linv \ff{\bu}{\bu
a}{\bu}$
do not vanish on $H$.
The $q$-Lie algebra, in virtue of {\sl Theorem 6.2}, is a $q$-Lie
subalgebra
of  $SO_{q,r=1}(N+2)$ .
It follows that the $\cchi{c_1}{c_2}$, $\cchi{b_1}{b_2}$
$q$-commutations
read as in eq. (\ref{Lierone}) with lower case indices:  they
give the $SO_{q,r=1}(N)$ $q$-Lie algebra. The remaining
commutations are [see (\ref{Lierone})]:
\eqa
& &\cchi{c_1}{c_2} \chi_{b_2} - {q_{c_1\bu} \over q_{c_2\bu}}
q_{b_2c_1} q_{c_2b_2} \chi_{b_2} \cchi{c_1}{c_2}=
{q_{c_1\bu} \over q_{c_2\bu}}[C_{b_2c_2} \chi_{c_1'}-\de^{c_1}_{b_2}
q_{c_2c_1} \chi_{c_2}]~, \label{Commchi0}\\
& &\chi_{c_2} \chi_{b_2} - {q_{b_2\bu} \over q_{c_2\bu}}
q_{c_2b_2} \chi_{b_2} \chi_{c_2}
=0~,\label{Commchi}\\
& &\cchi{c_1}{c_2} \cchi{\bu}{\bu} -\cchi{\bu}{\bu}
\cchi{c_1}{c_2}=0~~~~,~~~~
\chi_{c_2}\cchi{\bu}{\bu}-\cchi{\bu}{\bu}\chi_{c_2}=-\chi_{c_2}~
\ena
\noi where we have defined
$\chi_a \equiv  \cchi{\bu}{a}~.$
The exterior differential reads, $\forall a\in \ISO$
\eq
da=(\cchi{a}{b}*a)\Omega_{a}{}^b +
(\cchi{\bu}{b}*a)\Omega_{\bu}{}^b  +
(\cchi{\bu}{\bu}*a)\Omega_{\bu}{}^{\bu}~~;~~\mbox{
\sma{$a+b >  N+1$}}~
\en
where $\Omega_{a}{}^{b},~\Omega_{\bu}{}^{b},$
and $\Omega_{\bu}{}^{\bu}$ are the one-forms dual
to the tangent vectors (\ref{ISOtang1}) and (\ref{ISOtang2}).
Notice that the tangent vectors $\cchi{a}{b}$ and $\chi_{b}$
close on the $q$-Lie algebra (\ref{Commchi0}), (\ref{Commchi}) and
(\ref{Lierone}) with lower case indices.
This suggests
a reduction of the bicovariant calculus containing only the
$\cchi{a}{b}$ and $\cchi{\bu}{b}$ tangent vectors.
An explicit formulation, in agreement with \cite{Cas2}, is given in
\cite{altroarticolo}.
\sk
\noi{\bf Acknowledgments}
\sk
\noi {\small {$~~~$}The first author (P. A.) is  supported  by a  
joint
                fellowship
                University of California --  Scuola Normale  
Superiore,
                                                           Pisa,  
Italy and
by
Fondazione {\sl Angelo Della Riccia}, Firenze, Italy.
Part of this work has been accomplished through the Director,
Office of Energy Research,
Office af High Energy  and Nuclear Physics, Division of High Energy
Physics
of the U.S. Department of Energy under Contract DE-AC03-76SF00098.

 Work supported in part by EEC under TMR contract FMRX-CT96-0045.
}


\end{document}